# A GENERAL ASYMPTOTIC SCHEME FOR INFERENCE UNDER ORDER RESTRICTIONS

By D. Anevski and O. Hössjer

*Göteborg University and Stockholm University*

Limit distributions for the greatest convex minorant and its derivative are considered for a general class of stochastic processes including partial sum processes and empirical processes, for independent, weakly dependent and long range dependent data. The results are applied to isotonic regression, isotonic regression after kernel smoothing, estimation of convex regression functions, and estimation of monotone and convex density functions. Various pointwise limit distributions are obtained, and the rate of convergence depends on the self similarity properties and on the rate of convergence of the processes considered.

**1. Introduction.** Let $\{x_n\}_{n \geq 1}$ be a sequence of stochastic processes defined on an interval $J \subset \mathbb{R}$, and (a.s.) bounded from below on $J$. In this paper we consider the asymptotic behavior as $n \to \infty$ of the greatest convex minorant of $x_n$,

$$T_J(x_n) = \sup\{z; z : J \mapsto \mathbb{R}, z \text{ convex and } z \leq x_n\}, \tag{1}$$

at an interior point $t_0$ of $J$, as well as its derivative,

$$T_J(x_n)'(t) = \max_{v \leq t} \min_{u \geq t} \frac{x_n(u) - x_n(v)}{u - v}; \tag{2}$$

see Robertson, Wright and Dykstra [40]. Note that we use the convention $T_J(x)'(t) = T_J(x)'(t+)$ for any process $x$. The Pool Adjacent Violators Algorithm (PAVA) used to calculate $T$ can be found, for example, in [40].

The class of processes $x_n$ we consider includes partial sum and empirical processes for independent, weakly dependent and long range dependent data. The estimators (1) and (2) have several important applications, for instance, nonparametric regression and density estimation under order restrictions. The regression model has data $(y_i, t_i)_{i=1}^n$, satisfying

$$y_i = m(t_i) + \varepsilon_i,$$

where $m$ is the unknown regression function, $t_i = i/n$ are equidistant and $\{\varepsilon_i\}$ are error terms. If we restrict $m$ to be an increasing function, it is well known (cf. [10]) that the isotonic regression estimator

(3) $$\hat{m} = \arg\min\left\{\sum_{i=1}^n (y_i - z(t_i))^2 : z \text{ increasing}\right\}$$

is given by (2) at the observation points $\{t_i\}$, with $x_n$ the partial sum process formed by data. For independent and identically distributed errors $\{\varepsilon_i\}$, the asymptotic properties of $\hat{m}$ have been derived in [11, 47] and [33]. For instance, it follows from [11] that

(4) $$Cn^{1/3}(\hat{m}(t_0) - m(t_0)) \xrightarrow{\mathcal{L}} T(s^2 + B(s))'(0),$$

as $n \to \infty$, where $T = T_{\mathbb{R}}$, $B$ is a standard two-sided Brownian motion and $C$ depends on $m'(t_0)$ and $\sigma^2 = \text{Var}(\varepsilon_i)$. The right-hand side of (4) can also be replaced by $2\arg\min_{s \in \mathbb{R}}(s^2 + B(s))$, where we use the convention that, for any process $x$, $\arg\min_{s \in \mathbb{R}}(x(s))$ means the infimum of all points at which the minimum is attained.

In density estimation, data consists of a stationary process $\{t_i\}_{i=1}^n$ with an unknown marginal density function $f$. If $f$ is increasing and supported on a (finite or half-infinite) interval $J$, the nonparametric maximum likelihood estimate (NPMLE)

(5) $$\hat{f} = \arg\max\left\{\prod_{i=1}^n z(t_i) : z \text{ increasing}, z \geq 0 \text{ and } \int_J z(u)\, du = 1\right\}$$

for independent data can be written as $\hat{f} = T_J(x_n)'$, where $x_n$ is the empirical distribution; see [20]. Asymptotic properties of $\hat{f}$ have been obtained in [39] and [21]; see also [46]. In particular, (4) holds with $\hat{f}(t_0)$ and $f(t_0)$ in place of $\hat{m}(t_0)$ and $m(t_0)$, with $C$ a constant depending on $f'(t_0)$ and $f(t_0)$. Note also that increasing density estimation is related to unimodal density estimation; see [7] and references therein.

We propose to use $T_J(x_n)'$ as an estimator of $m$ and of $f$, also for dependent data; for the regression problem, $T_J(x_n)'$ minimizes the sum of squares in (3) no matter what dependence structure we have for $\{\varepsilon_i\}$, while the likelihood function is much more difficult to write down for dependent data; for density estimation the interpretation of $T_J(x_n)'$ is

$$T_J(x_n)' = \arg\min\left\{\sum_{i=1}^n (\tilde{x}_i - z(t_i))^2 w_i : z \text{ increasing}\right\},$$

where $\tilde{x}_i = (x_n(t_i) - x_n(t_{i-1}))/(t_i - t_{i-1})$ and $w_i = t_i - t_{i-1}$. Thus, $T_J(x_n)'$ is the weighted $\mathbf{l}^2$-projection of $(\tilde{x}_1, \ldots, \tilde{x}_n)$ on the convex set of increasing functions; see [40].

We review these results and show that the same limits are attained if data are weakly dependent and mixing. For long range dependent subordinated Gaussian data, in the regression problem we obtain a result reminiscent of (4), but with a different (nonpolynomial) convergence rate and with $B(\cdot)$ replaced with a process belonging to a class of long range dependent processes, which includes fractional Brownian motion and the Rosenblatt process; see [17]. In density estimation, $B(\cdot)$ is replaced by a straight line $Z \cdot s$ in (4) with $Z \sim N(0,1)$, in the cases for which we are able to check all the conditions. But since $T(s^2 + Zs)'(0) = (s^2 + Zs)'(0) = Z$, $\hat{f}$ is asymptotically normal in this case.

In [34] it is proposed, as an alternative to doing isotonic regression, to first smooth the data and then do isotonic regression, and the limit distribution is derived when using a kernel estimator with bandwidth $h \sim n^{-1/5}$ as smoother. We review these results, as well as state results for mixing and long range dependent data; however, we treat all possible choices of bandwidths $h$. An analogous approach is possible for density estimation; we, however, refrain from stating these results since it will be clear from the regression arguments how to proceed.

When estimating convex regression functions and density functions, the natural approaches would be to do convex regression or NPMLE of a convex density, respectively. An algorithm for convex regression has been proposed in [27], and a conjecture on the limit distribution can be found in [35]. In [30] an iterative algorithm for the NPMLE of a convex density and a conjecture on the limit distributions have been proposed; see also [2]. Finally, in [23, 24] the limit distributions for the convex regression and for the NPMLE of a convex density were derived.

As an alternative we propose the estimator $T_J(x_n)/c(x_n)$, where $x_n$ is a kernel estimate of either $m$ or $f$, and $c(x_n) = \int_J T_J(x_n)(u)\,du / \int_J x_n(u)\,du$. Thus, we obtain a convex function with the same integral over $J$ as $x_n$. The advantage over the regression and NMPLE approach is twofold: the PAVA algorithm used to calculate $T$ is noniterative and always converges, and in this paper we state the limit distributions of $T_J(x_n)$, both for the regression problem and for the density estimation problem, for weakly dependent data and long range dependent subordinated Gaussian data. The interpretation of $T_J(x_n)$ is the following: If $x_n'$ is piecewise continuous, then

$$(6) \qquad T_J(x_n)' = \arg\min\left\{\int_J (x_n'(u) - z(u))^2\,du : z \text{ increasing}\right\}$$

and, thus, $T_J(x_n)'$ is the $\mathbf{L}^2$-projection of $x_n'$ on the convex set of monotone functions; $T_J(x_n)$ is the primitive function of the solution to (6).

Our general convergence results can be written as

(7) $$d_n^{-p}(T_J(x_n)(t_0) - x_n(t_0)) \xrightarrow{\mathcal{L}} T(|s|^p + \tilde{v}(s))(0),$$

for the convex minorant of $x_n$, and

(8) $$d_n^{-p+1}(T_J(x_n)'(t_0) - x_{b,n}'(t_0)) \xrightarrow{\mathcal{L}} T(|s|^p + \tilde{v}(s))'(0),$$

for its derivative. Here $1 < p < \infty$ is a fixed number, $\tilde{v}$ a stochastic process reflecting the local behavior of $x_n$ around $t_0$ and $x_{b,n}$ is the deterministic part of $x_n$, for example, $E(x_n)$. The sequence $d_n \downarrow 0$ and $p$ determine the rate of convergence in (7) and (8). Values of $p$ different from 2 have previously been considered by Wright [47] and Leurgans [33], and arise, for example, in nonparametric regression when

$$m(t) - m(t_0) = a \operatorname{sgn}(t - t_0)|t - t_0|^{p-1} + o(|t - t_0|^{p-1}),$$

as $t \to t_0$ for some constant $a \neq 0$. The rate at which $d_n \downarrow 0$ depends on the rate of convergence of $x_n$ toward $x_{b,n}$ and on the local self similarity properties of $x_n$ around $t_0$.

Prakasa Rao [39] was the first to establish limit distributions for $T(x_n)'$ (with $x_n = F_n$ the empirical distribution function) and the approach presented in that paper has served as a model for later authors, first considering least convex minorants along a sequence of decreasing "truncated" intervals around $t_0$ and then establishing a truncation result saying that asymptotically the truncated intervals may replace $J$. Brunk [11] proved results for $T(x_n)'$ with $x_n$ the partial sum process, using similar techniques and relying on Prakasa Rao's result for the truncation reasoning. Wright [47] extended Brunk's result to cover monotone densities satisfying other smoothness assumptions, using a slightly different approach for the truncation proof. The methods used in these papers rely heavily on the fact that data are independent, using martingale results, and also on the fact that the limit process is a Brownian motion.

Leurgans [33] extended Wright's result to dependent data. The limit process is still assumed to be a Brownian motion, which could imply applications to weakly dependent data. However, the two applications given in [33] both deal with independent data (isotonic regression for independent and not identically distributed data, and isotonized quantile estimation for independent data). Next, Groeneboom [21] gave a different proof of Prakasa Rao's result, introducing strong approximation techniques (cf. [32]), and proved that the right-hand side of (3) is $2 \arg\min_{s \in \mathbb{R}}(s^2 + B(s))$; for the truncation result, a reference was made to Prakasa Rao's paper. Mammen [34] showed that a kernel estimate of $m$ with bandwidth $h \sim n^{-1/5}$ is first-order asymptotically equivalent to the estimate obtained by doing isotonic regression on the kernel estimate, thus obtaining the limit distribution for the isotonized kernel estimate. Wang [46] also used strong approximation to derive the limit distribution of the primitive function of the Grenander estimator.

A first example of a more general asymptotic theory of derivatives of least convex minorants can be found in [33], which potentially covers weakly dependent data. In our paper we treat both the convex minorant (1) and its derivative (2), arbitrary (nonpolynomial) sequences $d_n \downarrow 0$, as well as a large class of limit processes $\tilde{v}(\cdot)$, that is, not only Brownian motion. Thus, we are able to apply our general results also to estimation for dependent data (both short range and long range) and using estimates $x_n$ other than the partial sum process or empirical process, such as, for example, kernel estimates. Our method of proof is similar to the classical proof of Prakasa Rao [39], based on first considering least convex minorants along a sequence of decreasing "truncated" intervals around $t_0$ and then establishing a truncation result saying that asymptotically the truncated intervals may replace $J$. However, we decompose $x_n$ into a sum of a deterministic convex function $x_{b,n}$ and a stochastic part $v_n$. In this way we get very explicit regularity conditions that are possible to verify in a number of applications. Further, we use only weak convergence of a rescaled version of $v_n$ and do not refer to strong approximations, thereby obtaining greater generality. This relies on the application of the continuous mapping theorem and, thus, the continuity of the map $T_J : D(J) \mapsto C(J)$ is essential. Furthermore, we state conditions under which the continuous mapping theorem can be applied to the functional $x \mapsto T_J(x)'(t_0)$ (cf. Proposition 2). Such a condition automatically holds for Brownian motion and seems to have been implicitly assumed in previous work.

The article is organized as follows: Section 2 establishes the main convergence results (7) and (8) in Theorems 1 and 2, respectively. These results are then applied in Sections 3 and 4 to regression and density function estimation, respectively. In Section 5 a general formula is presented, which describes how $d_n$ depends on various properties of $x_n$, for example, local self similarity around $t_0$. In Section 6 we discuss possible extensions and generalizations. Finally, we have collected the proofs of the results in Section 2 and some technical empirical process and partial sum process results in the Appendix.

**2. Limit distributions.** Let $J \subset \mathbb{R}$ be a finite or infinite interval in $\mathbb{R}$ and define $D(J)$ as the space of functions $J \mapsto \mathbb{R}$ which are right continuous with left-hand limits.

Assume $\{x_n\}_{n \geq 1}$ is a sequence of stochastic processes on $D(J)$ for which we can write

(9) $$x_n(t) = x_{b,n}(t) + v_n(t), \qquad t \in J,$$

where $x_{b,n}$ is deterministic with $v_n$ also a member of $D(J)$. In this section we will derive limit distributions of $T_J(x_n)$ and $T_J(x_n)'$ for a large class of stochastic processes $x_n$. Our main assumptions on $x_n$ are that the process part $v_n$ can be rescaled in a way close to the self similarity property, and that the rescaled process converges weakly to some limit process. Given a sequence $d_n \downarrow 0$, we rescale $v_n$ locally around an interior point $t_0$ of $J$ according to

$$\tilde{v}_n(s;t_0) = d_n^{-p}(v_n(t_0 + sd_n) - v_n(t_0)),$$

where $1 < p < \infty$ is a fixed constant and $s \in J_{n,t_0} = d_n^{-1}(J - t_0)$. Thus, $\tilde{v}_n(\cdot;t_0) \in D(J_{n,t_0})$.

Many of the results on weak convergence are stated as results in $D[0,1]$ equipped with the Skorokhod metric. There are two reasons why this will not be appropriate for our needs. The first is that processes treated in our applications are not random elements of $D[0,1]$. For instance, $\tilde{v}_n(s;t_0)$ is defined on $D(J_{n,t_0}) = D[-a_n, b_n]$, where $a_n, b_n \to \infty$ as $n \to \infty$. The second reason is that the Skorokhod metric is too weak for the application we have in mind: the greatest convex minorant function $T: D[0,1] \mapsto C[0,1]$ will not be continuous if $D[0,1]$ is equipped with the Skorokhod topology. Thus, we would not be able to use the continuous mapping theorem to show limit distribution results for $T$ applied to a drift term plus a rescaled process.

The first problem is solved by working in $D(-\infty, \infty)$. For instance, $\tilde{v}_n(s)$ can be extrapolated according to

$$\tilde{v}_n(s;t_0) = \begin{cases} \tilde{v}_n(b_n+;t_0), & \text{if } s \geq b_n, \\ \tilde{v}_n(-a_n;t_0), & \text{if } s \leq -a_n. \end{cases}$$

Thus, $\tilde{v}_n(s;t_0)$ will lie in $D(-\infty, \infty)$ for all $n$. To deal with the second problem, we define a metric on $D(J)$ as follows: for $x, y \in D(J)$,

$$(10) \qquad \rho(x,y) = \sum_{k=1}^{\infty} 2^{-k} \frac{\rho_k(x,y)}{1 + \rho_k(x,y)},$$

where $\rho_k(x,y) = \sup_{s \in [-k,k] \cap J} |x(s) - y(s)|$, that is, we write $x_n \to x$ in $D(J)$ if, for each fixed $k$, $\sup_{[-k,k] \cap J} |x_n(s) - x(s)| \to 0$. Note that if $|J| < \infty$, then $\rho$ is equivalent to $\rho_J(x,y) = \sup_{s \in J} |x(s) - y(s)|$. By Theorem 23 in [38], page 108, weak convergence in $D(-\infty, \infty)$ is equivalent to weak convergence in $D[-k,k]$ of the processes restricted to $[-k,k]$, for every fixed $k$, where each $D[-k,k]$ of course is equipped with the sup-norm metric over $[-k,k]$. Note that with this metric the empirical process is not a measurable map if we use the Borel $\sigma$-algebra on $D[-k,k]$. If we instead use the $\sigma$-algebra generated by the open balls, the empirical process becomes measurable, and that assumption is also made in [38]. In that case, however, the continuous mapping theorem becomes somewhat more complicated, in that the set on which the function has all its continuity points should satisfy a certain regularity condition, as well as the usual demand that it have probability mass one. In the case of the functional $x \mapsto T(x)(t)$ that is not a problem, since this map is continuous everywhere; see (76) in Lemma A.1 in the sequel. However, in the case of the functional $x \mapsto T(x)'(t)$, it does pose a potential problem; see the proof of Proposition 2 and Note 2 in the sequel.

The next two assumptions are related to a local limit distribution result; see Lemma A.2 in Appendix A and the proof of Theorem 2.

ASSUMPTION A1 (Weak convergence of rescaled stochastic term). Assume there exists a stochastic process $\tilde{v}(\cdot;t_0) \neq 0$ such that

$$\tilde{v}_n(s;t_0) \xrightarrow{\mathcal{L}} \tilde{v}(s;t_0)$$

on $D(-\infty,\infty)$ as $n \to \infty$.

ASSUMPTION A2 (Bias term). Assume the functions $\{x_{b,n}\}_{n \geq 1}$ are convex. Put

(11)
$$g_n(s) = d_n^{-p}(x_{b,n}(t_0 + sd_n) - l_n(s)),$$
$$l_n(s) = x_{b,n}(t_0) + x'_{b,n}(t_0)sd_n,$$

for $s \in J_{n,t_0}$. Assume there is a constant $A > 0$ such that for each $c > 0$,

(12) $$\sup_{|s| \leq c} |g_n(s) - A|s|^p| \to 0,$$

as $n \to \infty$.

In applications we typically have a convex function $x_b$, such that either $x_{b,n} = x_b$ or $x_{b,n} \to x_b$ as $n \to \infty$, satisfying

$$x_b(t) = x_b(t_0) + x'_b(t_0)(t - t_0) + A|t - t_0|^p + o(|t - t_0|^p),$$

as $t \to t_0$. In particular, $A = \frac{1}{2}x''_b(t_0)$ if $p = 2$.

Define the rescaled function

(13) $$y_n(s) = g_n(s) + \tilde{v}_n(s;t_0).$$

The next two assumptions are related to a truncation result; see Lemma A.3 and Theorem A.1 in Appendix A.

ASSUMPTION A3 (Lower bound). For every $\delta > 0$, there are finite $0 < \tau = \tau(\delta)$ and $0 < \kappa = \kappa(\delta)$ such that

$$\liminf_{n \to \infty} P\left(\inf_{|s| \geq \tau}(y_n(s) - \kappa|s|) > 0\right) > 1 - \delta.$$

ASSUMPTION A4 (Small downdippings). Given $\varepsilon, \delta, \tilde{\tau} > 0$,

$$\limsup_{n \to \infty} P\left(\inf_{\tilde{\tau} \leq s \leq c} \frac{y_n(s)}{s} - \inf_{\tilde{\tau} \leq s} \frac{y_n(s)}{s} > \varepsilon\right) < \delta,$$

$$\limsup_{n \to \infty} P\left(\inf_{-c \leq s \leq -\tilde{\tau}} \frac{y_n(s)}{s} - \inf_{s \leq -\tilde{\tau}} \frac{y_n(s)}{s} < -\varepsilon\right) < \delta,$$

for all large enough $c > 0$.

We will now present a slightly less general but more transparent version of Assumptions A3 and A4, since in many of the applications it is possible to establish a separate restriction on the process part of $y_n$.

PROPOSITION 1. *Suppose Assumption A2 holds and that, for each $\varepsilon, \delta > 0$, there is a finite $\tau = \tau(\varepsilon, \delta)$ such that*

$$\limsup_{n \to \infty} P\left(\sup_{|s| \geq \tau} \left|\frac{\tilde{v}_n(s)}{g_n(s)}\right| > \varepsilon\right) < \delta. \tag{14}$$

*Then Assumptions A3 and A4 hold.*

Also the following assumption is related to the truncation results Lemma A.3 and Theorem A.1 in Appendix A.

ASSUMPTION A5 (Tail behavior of limit process). For each $\varepsilon, \delta > 0$, there is a $\tau = \tau(\varepsilon, \delta) > 0$ so that

$$P\left(\sup_{|s| \geq \tau} \left|\frac{\tilde{v}(s; t_0)}{A|s|^p}\right| > \varepsilon\right) \leq \delta.$$

THEOREM 1. *Let $t_0$ be fixed and suppose Assumptions A1, A2, A3, A4 and A5 hold. Then*

$$d_n^{-p}[T_J(x_n)(t_0) - x_n(t_0)] \xrightarrow{\mathcal{L}} T[A|s|^p + \tilde{v}(s; t_0)](0), \tag{15}$$

*with $A > 0$ as in Assumption A2, as $n \to \infty$.*

PROOF. Denote $T_c = T_{[-c,c]}, T = T_{\mathbb{R}}$ and $T_{c,n} = T_{[t_0 - cd_n, t_0 + cd_n]}$. Clearly,

$$d_n^{-p}(T_J(x_n)(t_0) - x_n(t_0)) = d_n^{-p}(T_J(x_n)(t_0) - T_{c,n}(x_n)(t_0))$$
$$+ d_n^{-p}(T_{c,n}(x_n)(t_0) - x_n(t_0)).$$

The truncation result in Lemma A.3 implies

$$d_n^{-p}(T_{c,n}(x_n)(t_0) - T_J(x_n)(t_0)) \xrightarrow{P} 0$$

if we first let $n \to \infty$ and then let $c \to \infty$. The local limit distribution result of Lemma A.2 implies that

$$d_n^{-p}(T_{c,n}(x_n)(t_0) - x_n(t_0)) \xrightarrow{\mathcal{L}} T_c[y(s)](0)$$

as $n \to \infty$, where

(16) $$y(s) = A|s|^p + \tilde{v}(s; t_0).$$

Then use Theorem A.1, Proposition 1 and Assumption A5 with $y_n(s) = y(s)$ to deduce

$$T_c(y(s))(0) - T(y(s))(0) \xrightarrow{P} 0$$

as $c \to \infty$. An application of Slutsky's theorem completes the proof. $\square$

Next we will study the limit distribution of the derivative $T(x_n)'$. There are some extra difficulties in this case. One is that the processes $x_n$ need not be differentiable. We therefore study the difference between $T(x_n)'$ and $x'_{b,n}$ directly.

Since the functional

(17) $$h : D[-c, c] \ni x \mapsto T(x)'(0)$$

is not continuous, the next assumption is essential.

ASSUMPTION A6. Suppose $y_n, y$ are defined in (13) and (16). Then

$$T_c(y_n)'(0) \xrightarrow{\mathcal{L}} T_c(y)'(0),$$

as $n \to \infty$, for each $c > 0$.

We need some simple condition in order to check Assumption A6.

PROPOSITION 2. *Assume $y$ takes its values in a separable set of completely regular points (cf. [38]), with probability one. Suppose Assumptions A1 and A2 hold and that for each $a \in \mathbb{R}$ and $c, \varepsilon > 0$,*

(18) $$P(y(s) - y(0) - \ as \ \geq \varepsilon|s| \ for \ all \ s \in [-c, c]) = 0.$$

*Then Assumption A6 holds.*

NOTE 1. Since $y(s) = \tilde{v}(s; t_0) + A|s|^p$ and $(A|s|^p)'(0) = 0$, (18) follows if we can prove

(19) $$P(\tilde{v}(\cdot; t_0) \in \Omega_c(a, \varepsilon)) = 0,$$

for each $a \in \mathbb{R}$ and $c, \varepsilon > 0$, with $\Omega_c(a, \varepsilon)$ defined in the proof of Proposition 2 in Appendix A. But (19) follows if we can find a random variable $Z$ (which may be deterministic) such that

$$P\left(\liminf_{s \to 0+} \frac{\tilde{v}(s; t_0) - Zs}{s} \leq 0\right) = 1, \tag{20}$$

$$P\left(\liminf_{s \to 0-} \frac{\tilde{v}(s; t_0) - Zs}{|s|} \leq 0\right) = 1. \tag{21}$$

Note that (20) and (21) hold if $\tilde{v}(s; t_0)$ is differentiable at 0 [take $Z = \tilde{v}'(0; t_0)$]. We can also make use of (the lower half of) the iterated logarithm law. Thus, with $Z = 0$, (20) and (21) follow if we can find a function $\psi : \mathbb{R} \setminus \{0\} \mapsto (0, \infty)$ such that

$$P\left(\liminf_{s \to 0+} \frac{\tilde{v}(s; t_0)}{\psi(s)} = -1\right) = 1,$$

$$P\left(\liminf_{s \to 0-} \frac{\tilde{v}(s; t_0)}{\psi(s)} = -1\right) = 1.$$

NOTE 2. If $y$ is continuous almost surely, the separability and complete regularity assumptions in Proposition 2 are satisfied; see Chapters 4 and 5 of [38]. All limit processes in this paper are almost surely continuous.

THEOREM 2. *Assume that Assumptions* A1–A6 *hold. Then*

$$d_n^{-p+1}[T(x_n)'(t_0) - x'_{b,n}(t_0)] \xrightarrow{\mathcal{L}} T(A|s|^p + \tilde{v}(s; t_0))'(0)$$

*as* $n \to \infty$. *Further, if*

$$P(T(A|s|^p + \tilde{v}(s; t_0))'(0) = a) = 0, \tag{22}$$

*then*

$$\lim_{n \to \infty} P\{d_n^{-p+1}[T(x_n)'(t_0) - x'_{b,n}(t_0)] < a\}$$

$$= P\left\{\underset{s \in \mathbb{R}}{\arg\min}(A|s|^p + \tilde{v}(s; t_0) - as) > 0\right\},$$

*with $A$ as in Assumption* A2.

PROOF. We start by proving a local limit distribution result. A $t$ varying in $I_n = [t_0 - cd_n, t_0 + cd_n]$ can be written as $t = t_0 + sd_n$ with $s \in [-c, c]$. Then

$$x_n(t_0 + sd_n) = v_n(t_0) + l_n(s) + d_n^p(g_n(s) + \tilde{v}_n(s; t_0)), \tag{23}$$

with $g_n, l_n$ defined in Assumption A2. We use the representation (23) and the chain rule to obtain

$$T_{c,n}(x_n)'(t_0) = x'_{b,n}(t_0) + d_n^{p-1} T_c(g_n(s) + \tilde{v}_n(s; t_0))'(0).$$

That is, Assumption A6 implies

$$d_n^{-p+1}(T_{c,n}(x_n)'(t_0) - x'_{b,n}(t_0)) = T_c(y_n)'(0) \xrightarrow{\mathcal{L}} T_c(y)'(0) \tag{24}$$

as $n \to \infty$, with $y_n, y$ defined in (13) and (16). Applying Lemma A.3 with $\Delta = 0$, we obtain

$$\lim_{c \to \infty} \limsup_{n \to \infty} P(d_n^{-p+1}|T_{c,n}(x_n)'(t_0) - T(x_n)'(t_0)| > \varepsilon) = 0. \tag{25}$$

Then, applying Theorem A.1 with $y_n(s) = y(s)$ and $I = \{0\}$, we get

$$\lim_{c \to \infty} P(|T_c(y(s))'(0) - T(y(s))'(0)| > \varepsilon) = 0. \tag{26}$$

Now (24), (25) and (26) and Slutsky's theorem prove the first part of the theorem.

For the second part of the theorem, we notice that if $P(T(y)'(0) = a) = 0$, then

$$\lim_{n \to \infty} P(d_n^{-p+1}(T(x_n)'(t_0) - x'_{b,n}(t_0)) < a) = P(T(y)'(0) < a).$$

Since $T(x)'$ is defined as the right-hand derivative of $T(x)$ and $\arg\min_{s \in \mathbb{R}}(x(s))$ as the infimum of all points at which the minimum is attained, it follows that

$$\{T(y)'(0) < a\} = \left\{ \arg\min_{s \in \mathbb{R}}(y(s) - as) > 0 \right\}. \tag{27}$$

By the first half of the theorem,

$$P(d_n^{-p+1}(T(x_n)'(t_0) - x'_{b,n}(t_0)) < a) \to P(T(y)'(0) < a),$$

if Assumption (22) holds, and this concludes the proof. □

In our applications the limit process will have stationary increments. Furthermore, it will be a two-sided version of a process defined on $\mathbb{R}^+$, and as such, its distribution will be unaffected by reflections in the $y$ axis through the origin. In these cases our results simplify.

ASSUMPTION A7 (Stationarity). The process $\tilde{v}(\cdot; t_0)$ has stationary increments, and

$$(\tilde{v}(s_1; t_0), \ldots, \tilde{v}(s_k; t_0)) \stackrel{\mathcal{L}}{=} (\tilde{v}(-s_1; t_0), \ldots, \tilde{v}(-s_k; t_0)),$$

for each $k$ and all $s_1, \ldots, s_k$.

COROLLARY 1. *Suppose Assumptions* A1–A7 *hold and $p = 2$. Then*

$$d_n^{-1}(T(x_n)'(t_0) - x'_{b,n}(t_0)) \xrightarrow{\mathcal{L}} 2\sqrt{A} \arg\min_{s \in \mathbb{R}} \left( s^2 + \tilde{v}\left(\frac{s}{\sqrt{A}}; t_0\right) \right)$$

*as $n \to \infty$, with $A$ as defined in Assumption* A2.

PROOF. We need to show that

$$\lim_{n\to\infty} P(d_n^{-1}(T(x_n)'(t_0) - x'_{b,n}(t_0)) < a)$$
$$= P\left(\arg\min_{s\in\mathbb{R}}\left(s^2 + \tilde{v}\left(\frac{s}{\sqrt{A}};t_0\right)\right) < \frac{a}{2\sqrt{A}}\right)$$

at each $a$ satisfying

(28) $$P\left(\arg\min_{s\in\mathbb{R}}\left(s^2 + \tilde{v}\left(\frac{s}{\sqrt{A}};t_0\right)\right) = \frac{a}{2\sqrt{A}}\right) = 0.$$

Note that

(29) $$\begin{aligned}P\left(\arg\min_{s\in\mathbb{R}}\left(s^2 + \tilde{v}\left(\frac{s}{\sqrt{A}};t_0\right)\right) < \frac{a}{2\sqrt{A}}\right)\\
= P\left(\arg\min_{s\in\mathbb{R}}\left(s^2 + \tilde{v}\left(\frac{s}{\sqrt{A}} + \frac{a}{\sqrt{A}};t_0\right)\right) > -\frac{a}{2\sqrt{A}}\right)\\
= P\left(\arg\min_{s\in\mathbb{R}}(y(s) - as) > 0\right),\end{aligned}$$

where the first equality follows by Assumption A7, and the second by a change of variables and completion of squares. Putting $h_a(y) = \mathbb{1}_{\{\arg\min(y(s)-as)>0\}}$, we can rewrite (28) as

(30) $$\lim_{\varepsilon\to 0} Eh_{a+\varepsilon}(y) = Eh_a(y).$$

Note also that $h_{a+\varepsilon}(y) \uparrow$ when $\varepsilon \downarrow 0$. Let $D = \{z: \lim_{\varepsilon\to 0} h_{a+\varepsilon}(z) \neq h_a(z)\}$. Then if $P(D) = 0$, we have

$$Eh_{a+\varepsilon}(y) = Eh_{a+\varepsilon}(y)\mathbb{1}_{\{D^c\}}(y) \uparrow Eh_a(y)\mathbb{1}_{\{D^c\}}(y) = Eh_a(y)$$

as $\varepsilon \downarrow 0$ by monotone convergence, and, thus, (30) holds. But (27) implies

$$h_{a+\varepsilon}(z) = 1 \iff T(z)'(0) < a + \varepsilon.$$

Thus

$$D = \{z: T(z)'(0) = a\},$$

and the latter part of Theorem 2 completes the proof. □

Let $S(x_n)$ denote the least concave majorant of $x_n$. Limit distribution results for $S(x_n)$ and $S(x_n)'$ now follow easily by noting that $S(x_n) = -T(-x_n)$.

In the next two sections we will consider various applications of Theorems 1 and 2 when $p = 2$. Applications for other $p > 1$ in the independent data case are treated in [47] and [33].

**3. Regression.** Assume $m$ is a function on the interval $J = [0,1] \subset \mathbb{R}$, and $(y_i, t_i), i = 1, \ldots, n$, are pairs of data satisfying

$$y_i = m(t_i) + \varepsilon_i, \tag{31}$$

where the $t_i = i/n$ are the design points, that is, we have an equispaced design. For later convenience, we define the error terms $\varepsilon_i$ for all integers, and assume that $\{\varepsilon_i\}_{i=-\infty}^{\infty}$ form a stationary sequence of random variables with $E(\varepsilon_i) = 0$ and $\text{Var}(\varepsilon_i) = \sigma^2 < \infty$. Let $\sigma_n^2 = \text{Var}(\sum_{i=1}^{n} \varepsilon_i)$. Then the two-sided partial sum process $w_n$ is defined by

$$w_n\left(t_i + \frac{1}{2n}\right) = \begin{cases} \dfrac{1}{\sigma_n}\left(\dfrac{\varepsilon_0}{2} + \sum_{j=1}^{i} \varepsilon_i\right), & i = 0, 1, 2, \ldots, \\ \dfrac{1}{\sigma_n}\left(-\dfrac{\varepsilon_0}{2} - \sum_{j=i+1}^{-1} \varepsilon_i\right), & i = -1, -2, \ldots, \end{cases}$$

and linearly interpolated between these points. This process is right continuous with left-hand limits, so it lies in the space $D(-\infty, \infty)$.

The dependence structure for the random parts, the $\varepsilon_i$, will determine the limit distribution. Let $\text{Cov}(k) = E(\varepsilon_1 \varepsilon_{1+k})$ denote the covariance function. Then it is possible to distinguish between three cases [of which (i) is a special case of (ii)]:

(i) Independence: the $\varepsilon_i$ are independent.
(ii) Weak dependence: $\sum_k |\text{Cov}(k)| < \infty$.
(iii) Strong (long range) dependence: $\sum_k |\text{Cov}(k)| = \infty$.

The first two cases are similar in the sense that, for these, $w_n$ has the same limit distribution, namely, the Brownian motion. For the case of long range dependence, the limit distributions are very different. Also, this case is the most awkward to work with, and limit distribution results are known only for subordinated processes, that is, when $\varepsilon_i$ is a function of an underlying process with a parametric law. We will treat only subordinated Gaussian processes when the underlying process is Gaussian. All results stated will be for processes in $D(-\infty, \infty)$ with the uniform metric on compacta defined in (10), and the $\sigma$-algebra generated by the open balls.

Most of the limit results stated for partial sum processes are results for processes in $D[0, 1]$ equipped with the Skorokhod metric. An examination of the proofs of the limit distribution results for $D[0, 1]$ shows that there is nothing special about $[0, 1]$; it can be replaced by $[0, k]$, for any finite $k$. This means that the results can be seen as results for $D[0, k]$, with the Skorokhod metric. If the limit process is in $C[0, k]$ a.s., we can use the Skorokhod–Dudley theorem to get new random processes converging almost surely, so in the Skorokhod metric on a set with probability one. But convergence in that metric toward a continuous function implies convergence in the supnorm-metric, and this implies weak convergence in $D[0,k]$ with the supnorm-topology. Finally, this is made into a result for $D[-k,k]$, for the two-sided partial sum process $w_n$.

When the $\varepsilon_i$ are independent, we have the classical Donsker theorem (cf. [8]), implying that

$$w_n \stackrel{\mathcal{L}}{\to} B, \tag{32}$$

as $n \to \infty$, with $B$ a two-sided standard Brownian motion on $D(-\infty, \infty)$.

Next we treat weakly dependent data. The notion of weak dependence can be formalized in several ways. We will use mixing conditions; for a survey see [9]. Define the $\sigma$-algebras

$$\mathcal{F}_k = \sigma\{\varepsilon_i : i \leq k\},$$
$$\bar{\mathcal{F}}_k = \sigma\{\varepsilon_i : i \geq k\},$$

where $\sigma\{\varepsilon_i : i \in I\}$ denotes the $\sigma$-algebra generated by $\{\varepsilon_i : i \in I\}$.

DEFINITION 1. The stationary sequence $\{\varepsilon_i\}$ is said to be $\phi$-mixing or $\alpha$-mixing, respectively, if there is a function $\phi(n)$ or $\alpha(n) \to 0$ as $n \to \infty$, such that

$$\sup_{A \in \bar{\mathcal{F}}_n} |P(A|\mathcal{F}_0) - P(A)| \leq \phi(n),$$
$$\sup_{A \in \mathcal{F}_0, B \in \bar{\mathcal{F}}_n} |P(AB) - P(A)P(B)| \leq \alpha(n).$$

Mixing conditions say that elements in sequences are almost independent if they are far away from each other. There are other ways to model weak dependence, such as the notion of mixingales introduced in [36], which is a special case of the processes treated in [26]. See also the results for short range dependent subordinated Gaussian sequences in [15].

Introduce

$$\kappa^2 = \text{Cov}(0) + 2\sum_{k=1}^{\infty} \text{Cov}(k) \tag{33}$$

whenever the limit exists. The following results for mixing sequences are adapted from [37] and [26].

ASSUMPTION A8 ($\phi$-mixing). Assume $\{\varepsilon_i\}_{i \in \mathbf{Z}}$ is a stationary $\phi$-mixing sequence with $E\varepsilon_i = 0$ and $E\varepsilon_i^2 < \infty$. Assume further $\sum_{k=1}^{\infty} \phi(k)^{1/2} < \infty$ and $\kappa > 0$ in (33).

Note that $\kappa^2$ exists and that

(34) $$\frac{\sigma_n^2}{n} \to \kappa^2,$$

as $n \to \infty$ by Lemmas 20.1 and 20.3 in [8], if Assumption A8 is satisfied. In [26] it is shown that Donsker's result (32) is implied by Assumption A8 and also by several other combinations of assumptions.

ASSUMPTION A9 ($\alpha$-mixing). Assume $\{\varepsilon_i\}_{i \in \mathbf{Z}}$ is a stationary $\alpha$-mixing sequence with $E\varepsilon_i = 0$ and $E\varepsilon_i^4 < \infty$, $\kappa > 0$ in (33) and $\sum_{k=1}^{\infty} \alpha(k)^{1/2-\varepsilon} < \infty$, for some $\varepsilon > 0$.

From Lemma 20.1 in [8] and Theorem 17.2.2 in [29] it follows that $\kappa^2$ exists and that (34) holds, if Assumption A9 is satisfied. The results of Peligrad [37] imply that if Assumption A9 holds, then Donsker's result (32) follows.

To treat long range dependent data, assume $\{\xi_i\}_{i \in \mathbf{Z}}$ is a stationary Gaussian process with mean zero and covariance function $\mathrm{Cov}(k) = E(\xi_i \xi_{i+k})$ such that $\mathrm{Cov}(0) = 1$ and $\mathrm{Cov}(k) = k^{-d} l_0(k)$, where $l_0$ is a function slowly varying at infinity and $0 < d < 1$ is fixed. For a review of long range dependence, see [6].

Let $g : \mathbb{R} \mapsto \mathbb{R}$ be a measurable function and define $\varepsilon_i = g(\xi_i)$. Then we can expand $g(\xi_i)$ in Hermite polynomials

$$g(\xi_1) = \sum_{k=r}^{\infty} \frac{1}{k!} \eta_k h_k(\xi_1),$$

with equality holding as a limit in $L^2(\phi)$, with $\phi$ the standard Gaussian density function. Here $h_k$ are the Hermite polynomials of order $k$, the functions

$$\eta_k = E(g(\xi_1) h_k(\xi_1)) = \int g(u) h_k(u) \phi(u)\, du,$$

are the $L^2(\phi)$-projections on $h_k$, and $r$ is the index of the first nonzero coefficient in the expansion. Assuming that $0 < dr < 1$, the sequence $\{\varepsilon_i\}$ also exhibits long range dependence. In this case we say the sequence $\{\varepsilon_i\}$ is subordinated Gaussian long range dependent with parameters $d$ and $r$.

The results of Taqqu [43, 44] show that

$$\sigma_n^{-1} \sum_{i \leq nt} g(\xi_i) \xrightarrow{\mathcal{L}} z_{r,\beta}(t)$$

in $D[0,1]$ equipped with the Skorokhod topology. Lemma 3.1 and Theorem 3.1 in [43] show that the variance is $\sigma_n^2 = \mathrm{Var}(\sum_{i=1}^n g(\xi_i)) = \eta_r^2 n^{2-rd} l_1(n)(1 + o(1))$, where

(35) $$l_1(k) = \frac{2}{r!(1-rd)(2-rd)} l_0(k)^r.$$

The limit process $z_{r,\beta}$ is in $C[0,1]$ a.s., and is self similar with parameter

$$\beta = 1 - rd/2. \tag{36}$$

That is, the processes $z_{r,\beta}(\delta t)$ and $\delta^\beta z_{r,\beta}(t)$ have the same finite-dimensional distributions for all $\delta > 0$.

The limit process can, for arbitrary $r$, be represented by Wiener–Itô–Dobrushin integrals as in [17]; see also the representation given in [44]. The process $z_{1,\beta}(t)$ is fractional Brownian motion, $z_{2,\beta}(t)$ is the Rosenblatt process, and the processes $z_{r,\beta}(t)$ are all non-Gaussian for $r \geq 2$; see Taqqu [43]. This implies that, under the above assumptions,

$$w_n \xrightarrow{\mathcal{L}} B_{r,\beta} \tag{37}$$

in $D(-\infty, \infty)$, as $n \to \infty$, where $B_{r,\beta}$ are the two-sided versions of the processes $z_{r,\beta}$.

3.1. *Isotonic regression.* Assume the regression function $m$ in (31) satisfies $m \in \mathcal{F} = \{$increasing functions$\}$. The problem of minimizing the sum of squares $\sum_{i=1}^n (y_i - m(t_i))^2$ over the class $\mathcal{F}$ is known as the isotonic regression problem. The nonparametric least squares estimator is obtained as

$$\hat{m} = T_{[0,1]}(x_n)'$$

(cf., e.g., [40]), where $x_n$ is defined as follows: Let $\tilde{n} = \tilde{n}(t) = \lfloor nt - 1/2 \rfloor$ and put

$$x_n(t) = n^{-1} \sum_{i=1}^{\tilde{n}} y_i + \frac{(nt - 1/2) - \tilde{n}}{n} y_{\tilde{n}+1}, \qquad t \in [0,1].$$

The limit distribution of $T_{[0,1]}(x_n)'$ is known in the case of independent data and is included in Theorem 3 in [35]; note also the results in [33] and [47].

Actually $T_{[0,1]}(x_n)'$ is the solution to the isotonic regression problem, no matter what the dependence structure is, and we will derive the limit distributions also in the weakly and long range dependent cases.

We can write

$$x_n(t) = x_{b,n}(t) + v_n(t),$$

with

$$x_{b,n}(t) = n^{-1} \sum_{i=1}^{\tilde{n}} m(t_i) + \frac{(nt - 1/2) - \tilde{n}}{n} m(t_{\tilde{n}+1}),$$

$$v_n(t) = n^{-1} \sum_{i=1}^{\tilde{n}} \varepsilon_i + \frac{(nt - 1/2) - \tilde{n}}{n} \varepsilon_{\tilde{n}+1}.$$

Note that $x_{b,n}$ is convex (recall that we assume $p = 2$ in Sections 3 and 4) and that, because of the stationarity of $\{\varepsilon_i\}_{i=-\infty}^{\infty}$,

$$\begin{aligned}\tilde{v}_n(s;t) &= d_n^{-2}(v_n(t+sd_n) - v_n(t)) \\ &= d_n^{-2} n^{-1} \sigma_{\hat{n}}(w_{\hat{n}}(td_n^{-1} + s) - w_{\hat{n}}(td_n^{-1})) \\ &\stackrel{\mathcal{L}}{=} d_n^{-2} n^{-1} \sigma_{\hat{n}} w_{\hat{n}}(s),\end{aligned}$$

where $\hat{n} = nd_n$ and the last equality in distribution holds exactly when $t = t_i$ for any $i$ and asymptotically for all $t$. Since we know that, under the appropriate assumptions, $w_{\hat{n}} \stackrel{\mathcal{L}}{\to} w$ in $D(-\infty, \infty)$ for some process $w$, we need to choose $d_n$ in such a way that $d_n^{-2} n^{-1} \sigma_{\hat{n}} \to c$ for some constant $0 < c < \infty$. Thus,

$$\tilde{v}_n(s;t_0) \stackrel{\mathcal{L}}{\to} cw(s) =: \tilde{v}(s;t_0)$$

in $D(-\infty, \infty)$.

THEOREM 3. *Assume $m$ is increasing with $m'(t_0) > 0$ and $t_0 \in (0,1)$. Let $\hat{m}(t_0) = T_{[0,1]}(x_n)'(t_0)$ be the solution to the isotonic regression problem. Suppose that one of the following conditions holds:*

(i) $\{\varepsilon_i\}$ *are independent and identically distributed with $E\varepsilon_i = 0$ and* $\text{Var}(\varepsilon_i) = \sigma^2 < \infty$;

(ii) *Assumption* A8 *or* A9 *holds, $\sigma_n^2 = \text{Var}(\sum_{i=1}^n \varepsilon_i)$ and define $\kappa^2$ as in* (33);

(iii) $\varepsilon_i = g(\xi_i)$ *is a long range dependent subordinated Gaussian sequence with parameters $d$ and $r$, and $\beta$ as in* (36).

*Then, correspondingly, we obtain*

$$d_n^{-1} c_1(t_0)(\hat{m}(t_0) - m(t_0)) \stackrel{\mathcal{L}}{\to} \underset{s \in \mathbb{R}}{\arg\min}(s^2 + \tilde{v}(s)),$$

$$d_n^{-2} c_2(t_0)\left(\int_0^{t_0} \hat{m}(s)\,ds - x_n(t_0)\right) \stackrel{\mathcal{L}}{\to} T(s^2 + \tilde{v}(s))(0),$$

*as $n \to \infty$ with, respectively:*

(i) $\tilde{v} = B, d_n = n^{-1/3}, c_1(t_0) = 2^{-2/3} m'(t_0)^{-1/3} \sigma^{-2/3}$, $c_2(t_0) = 2^{-1/3} \times m'(t_0)^{1/3} \sigma^{-4/3}$;

(ii) $\tilde{v} = B, d_n = n^{-1/3}, c_1(t_0) = 2^{-2/3} m'(t_0)^{-1/3} \kappa^{-2/3}$, $c_2(t_0) = 2^{-1/3} \times m'(t_0)^{1/3} \kappa^{-4/3}$;

(iii) $\tilde{v} = B_{r,\beta}, d_n = l_2(n) n^{-rd/(2+rd)}, c_1(t_0) = 2^{-1/(2-\beta)} m'(t_0)^{(\beta-1)/(2-\beta)} \times |\eta_r|^{-1/(2-\beta)}, c_2(t_0) = 2^{-\beta/(2-\beta)} m'(t_0)^{\beta/(2-\beta)} |\eta_r|^{-2/(2-\beta)}$;

and $l_2$ is a slowly varying function related to $l_1$ as shown in the proof below. Moreover,

$$n\sigma_{\tilde{n}}^{-1} \int_0^{t_0} (\hat{m}(s) - m(s))\, ds \xrightarrow{\mathcal{L}} w(1) \tag{38}$$

as $n \to \infty$, where $w(t) = B(t)$ in the cases (i) and (ii) and $w(t) = B_{r,\beta}(t)$ in the case (iii), and $\tilde{n} = \lfloor nt - 1/2 \rfloor$.

PROOF. (i) (The independent case) We have $\sigma_{\hat{n}}^2 = \sigma^2 \hat{n} = \sigma^2 n d_n$, which implies that we can choose $d_n = n^{-1/3}$, so that $c = d_n^{-2} n^{-1} \sigma_{\hat{n}} = \sigma$. The rescaled process is $y_n(s) = g_n(s) + \tilde{v}_n(s;t)$, where $\tilde{v}_n, g_n$ are defined in Section 2. From Donsker's theorem (32), it follows that $\tilde{v}_n(s;t_0) \xrightarrow{\mathcal{L}} \sigma B(s)$ as $n \to \infty$ on $D(-\infty, \infty)$ and, thus, Assumption A1 is satisfied. Next define $\tilde{m}_n(t) = m(t_i)$ when $t_i - 1/(2n) < t \le t_i + 1/(2n)$, so that $x_{b,n}(t) = \int_0^t \tilde{m}_n(u)\, du$. Then

$$g_n(s) = d_n^{-2} \int_{t_0}^{t_0+sd_n} (\tilde{m}_n(u) - \tilde{m}_n(t_0))\, du$$

$$= d_n^{-2} \int_{t_0}^{t_0+sd_n} (m(u) - m(t_0))\, du + r_n(s),$$

where the first term converges toward $As^2$ uniformly for $s$ on compacta, with $A = m'(t_0)/2$, and

$$\sup_{|s| \le c} |r_n(s)| \le 2cd_n^{-1} \sup_{|u-t_0| \le sd_n} |m(u) - m_n(u)| = O(n^{-1} d_n^{-1}) = o(1),$$

since $n^{-1} d_n^{-1} = c d_n \sigma_{\hat{n}}^{-1}(1 + o(1)) \to 0$, because $d_n \to 0$ and $\sigma_{\hat{n}} \to \infty$ and, thus, Assumption A2 holds. Assumptions A3 and A4 follow by Proposition 1 and Lemma B.1 in Appendix B. Assumptions A5, A6 and A7 hold by properties of the Brownian motion; see [42] for an LIL for Brownian motion which shows Assumption A6 via Proposition 2 and Note 1. Thus, from Theorem 1,

$$n^{2/3}\left(\int_0^{t_0} \hat{m}(s)\, ds - x_n(t_0)\right) \xrightarrow{\mathcal{L}} T(As^2 + \sigma B(s))(0)$$

$$\stackrel{\mathcal{L}}{=} A^{-1/3} \sigma^{4/3} T(s^2 + B(s))(0)$$

as $n \to \infty$, where the equality follows from the self similarity of Brownian motion. Furthermore, Corollary 1 implies

$$n^{1/3}(\hat{m}(t_0) - m(t_0)) \xrightarrow{\mathcal{L}} 2A^{1/2} \arg\min_{s \in \mathbb{R}}\left(s^2 + \sigma B\left(\frac{s}{\sqrt{A}}\right)\right)$$

$$\stackrel{\mathcal{L}}{=} 2A^{1/3} \sigma^{2/3} \arg\min_{s \in \mathbb{R}}(s^2 + B(s))$$

as $n \to \infty$, where the equality follows from the self similarity of Brownian motion, and the proof is complete.

(ii) (The mixing case) Choosing $\hat{n} = nd_n$, we get, as in the independent data case, $d_n^{-2} n^{-1} \sigma_{\hat{n}} \to \kappa$, so that $\tilde{v}(s; t_0) = \kappa B(s)$, with the choice $d_n = n^{-1/3}$. The rest of the proof goes through as for independent data.

(iii) (The long range dependent case) In this case $\sigma_{\hat{n}}^2 = \eta_r^2 (nd_n)^{2-rd} l_1(\hat{n})$. We choose $d_n$ as

$$|\eta_r| = d_n^{-2} n^{-1} \sigma_{\hat{n}} = d_n^{-2} n^{-1} |\eta_r| (nd_n)^{1-rd/2} l_1(nd_n)^{1/2}$$
$$\iff \quad d_n^{1+rd/2} = n^{-rd/2} l_1(nd_n)^{1/2}$$
$$\iff \quad d_n = n^{-rd/(2+rd)} l_2(n),$$

where $l_2$ is another function slowly varying at infinity. Thus,

$$\tilde{v}_n(s; t_0) = |\eta_r| w_{\hat{n}}(s) \xrightarrow{\mathcal{L}} \tilde{v}(s; t_0) = |\eta_r| B_{r,\beta},$$

on $D(-\infty, \infty)$, as $n$ (and $\hat{n}$) $\to \infty$, and Assumption A1 holds. Assumption A2 is proved as for independent data and Assumptions A3 and A4 follow from Proposition 1 and Lemma B.1. Also Assumptions A5, A6 and A7 follow from the properties of $B_{r,\beta}$; see Proposition 2 for Assumption A6. The assumptions of Theorem 1 are therefore satisfied and

$$d_n^{-2} \left( \int_0^{t_0} \hat{m}(s) \, ds - x_n(t_0) \right) \xrightarrow{\mathcal{L}} T(As^2 + |\eta_r| B_{r,\beta}(s))(0)$$
$$\stackrel{\mathcal{L}}{=} A^{-\beta/(2-\beta)} |\eta_r|^{2/(2-\beta)} T(s^2 + B_{r,\beta}(s))(0)$$

as $n \to \infty$, where the equality follows from the self similarity of $B_{r,\beta}$. Furthermore, Corollary 1 implies

$$d_n^{-1}(\hat{m}(t_0) - m(t_0)) \xrightarrow{\mathcal{L}} 2A^{1/2} \arg\min_{s \in \mathbb{R}} \left( s^2 + |\eta_r| B_{r,\beta}\left( \frac{s}{\sqrt{A}} \right) \right)$$
$$\stackrel{\mathcal{L}}{=} 2A^{(1-\beta)/(2-\beta)} |\eta_r|^{1/(2-\beta)} \arg\min_{s \in \mathbb{R}} (s^2 + B_{r,\beta}(s)),$$

where the equality follows from the self similarity of $B_{r,\beta}$.

To show (38), note that, with $M(t) = \int_0^t m(s) \, ds$,

$$x_n(t_0) - M(t_0) = v_n(t_0) + (x_{b,n}(t_0) - M(t_0))$$
$$= n^{-1} \sigma_{\tilde{n}} w_n(1) + O_p(n^{-1}),$$

implying that

$$n \sigma_{\tilde{n}}^{-1} (x_n(t_0) - M(t_0)) \xrightarrow{\mathcal{L}} w(1)$$

as $n \to \infty$. This proves the theorem. $\square$

3.2. *Estimating a convex regression function.* Assume the regression function $m$ in (31) belongs to the class $\mathcal{F}_2 = \{\text{convex functions}\}$. One natural way to estimate $m$ based on the data is to do convex regression, that is, to minimize the sum of squares over the class of convex functions. Algorithms for the convex regression problem are given in [27] and [35], and the limit distribution for independent data is presented in [23, 24]. We present here an estimator of a convex regression function for which we are able to give the limit distributions also in the weakly dependent and long range dependent cases.

Define $\bar{y}_n : [1/n, 1] \mapsto \mathbb{R}$ by linear interpolation of the points $\{(t_i, y_i)\}_{i=1}^n$, and let

$$(39) \qquad x_n(t) = h^{-1} \int k((t-u)/h) \bar{y}_n(u) \, du$$

be the Gasser–Müller kernel estimate of $m(t)$ (see [19]), where $k$ is a symmetric density in $\mathbf{L}^2(\mathbb{R})$ with compact support; for simplicity, take $\mathrm{supp}(k) = [-1, 1]$; $k$ is called the kernel function. Let $h$ be the bandwidth, for which we assume that $h \to 0, nh \to \infty$. The exact choice of $h$ will be affected by the dependence structure of $\{\varepsilon_i\}$.

To define a convex estimator of $m$, we put

$$(40) \qquad \tilde{m}(t) = \frac{T_{[0,1]}(x_n)(t)}{c(x_n)}, \qquad t \in J,$$

where $c(x_n) = \int_J T_{[0,1]}(x_n)(t) \, dt (\int_J x_n(t) \, dt)^{-1}$ is a normalization constant that ensures $\int_J \tilde{m}(t) \, dt = \int_J x_n(t) \, dt$. We will confine ourselves to studying the asymptotics of $T_{[0,1]}(x_n)$, that is, the behavior of $\tilde{m}$ before normalization. Kernel regression estimation for long range dependent errors is considered in [13, 14].

Clearly, $x_n(t) = x_{b,n}(t) + v_n(t)$, with

$$x_{b,n}(t) = h^{-1} \int k\left(\frac{t-u}{h}\right) \bar{m}_n(u) \, du,$$

$$v_n(t) = h^{-1} \int k\left(\frac{t-u}{h}\right) \bar{\varepsilon}_n(u) \, du,$$

where the functions $\bar{m}_n$ and $\bar{\varepsilon}_n$ are obtained by linear interpolation of $\{(t_i, m(t_i))\}_{i=1}^n$ and $\{(t_i, \varepsilon_i)\}_{i=1}^n$, respectively. For the deterministic term, $x_{b,n}(t) \to x_b(t) = m(t)$, as $n \to \infty$. Note that $\bar{m}_n$, and thus also $x_{b,n}$, is convex.

Put

$$(41) \qquad \bar{w}_n(t) = \frac{n}{\sigma_n} \int_0^t \bar{\varepsilon}_n(u) \, du.$$

Since $\text{supp}(k) = [-1, 1]$ and if $t \in (1/n + h, 1 - h)$, from a partial integration and change of variable, we obtain

$$v_n(t) = \frac{\sigma_n}{nh} \int k'(u) \bar{w}_n(t - uh) \, du.$$

It can be shown that $\bar{w}_n$ and $w_n$ are asymptotically equivalent for all dependence structures treated in this paper. This will henceforth be tacitly assumed.

Recall that for the rescaling of $v_n$ we need to choose $d_n$ in a correct way. Having done that choice, depending on the relation between the rate of convergence to zero of the bandwidth and of $d_n$, we get different limit results for $T(x_n)$. We have three subcases: $d_n = h, d_n/h \to 0$, or $d_n/h \to \infty$ as $n \to \infty$.

3.2.1. *The case $d_n = h$.* For $s > 0$, we rescale as

$$\tilde{v}_n(s;t) = d_n^{-2}(nh)^{-1}\sigma_{\hat{n}} \int (\bar{w}_{\hat{n}}(h^{-1}t + s - u)$$

$$- \bar{w}_{\hat{n}}(h^{-1}t - u))k'(u) \, du$$

$$\stackrel{\mathcal{L}}{=} d_n^{-2}(nh)^{-1}\sigma_{\hat{n}} \int (\bar{w}_{\hat{n}}(s - u) - \bar{w}_{\hat{n}}(-u))k'(u) \, du,$$

with $\hat{n} = nh$, where the last equality holds exactly only for $t = t_i$ and asymptotically otherwise. Note that the right-hand side holds also for $s < 0$.

Assume $d_n = h$ is such that

(42) $$d_n^{-2}(nh)^{-1}\sigma_{\hat{n}} = d_n^{-3}n^{-1}\sigma_{\hat{n}} \to c > 0.$$

Then, under conditions given in the beginning of this chapter, $w_n \stackrel{\mathcal{L}}{\to} w$ in $D(-\infty, \infty)$, using the supnorm over compacta metric. Note that if $k'$ is bounded and $k$ has compact support, the map

$$D(-\infty, \infty) \ni z(s) \mapsto \int (z(s - u) - z(-u))k'(u) \, du \in D(-\infty, \infty)$$

is continuous, using the supnorm over compacta. Thus, the continuous mapping theorem implies that

(43) $$\tilde{v}_n(s;t) \stackrel{\mathcal{L}}{\to} \tilde{v}(s;t) = c \int (w(s - u) - w(-u))k'(u) \, du.$$

Define $\hat{m} = T_{[0,1]}(x_n(t))$, and note in the following theorem that the rate $n^{-2/5}$ in the independent data case is the same as the rate in the limit distribution result for the convex regression; see [23, 24].

THEOREM 4. *Assume $m$ is convex with $m''(t_0) > 0$ and $t_0 \in (0,1)$. Let $x_n$ be the kernel estimate of $m$ defined in (39), with a nonnegative and compactly supported kernel $k$ such that $k'$ is bounded, and with bandwidth $h$ specified below. Suppose that one of the following conditions holds:*

  (i) *$\{\varepsilon_i\}$ are independent and identically distributed with $E\varepsilon_i = 0$ and $\sigma^2 = \operatorname{Var}(\varepsilon_i) < \infty$ and we choose $h = an^{-1/5}$, where $a > 0$ is an arbitrary constant;*
  (ii) *Assumption A8 or A9 holds, $\sigma_n^2 = \operatorname{Var}(\sum_{i=1}^n \varepsilon_i)$ and $\kappa^2$ is defined in (33), and we choose $h = an^{-1/5}$, where $a > 0$ is an arbitrary constant;*
  (iii) *$\varepsilon_i = g(\xi_i)$ is a long range dependent subordinated Gaussian sequence with parameters $d$ and $r$, and $\beta$ as in (36) and we choose $h = l_2(n;a)n^{-rd/(4+rd)}$, where $a > 0$ and $n \mapsto l_2(n;a)$ is a slowly varying function defined in the proof below.*

*Then, correspondingly, we obtain*

$$d_n^{-2}(T_{[0,1]}(x_n)(t_0) - m(t_0)) \xrightarrow{\mathcal{L}} \tfrac{1}{2}m''(t_0)\int u^2 k(u)\,du + c\int k'(u)w(-u)\,du$$
$$+ T(\tfrac{1}{2}m''(t_0)s^2 + \tilde{v}(s;t_0))(0),$$

*as $n \to \infty$, where $\tilde{v}(s;t)$ is defined in (43), $d_n = h$ and, respectively:*

  (i) $w = B, c = \sigma a^{-5/2}$,
  (ii) $w = B, c = \kappa a^{-5/2}$,
  (iii) $w = B_{r,\beta}, c = |\eta_r|a$.

PROOF. (i) (Independent case) We have $\sigma_{\hat{n}}^2 = \sigma^2 n d_n$. Thus, $d_n^{-2}(nh)^{-1}\sigma_{\hat{n}} = \sigma n^{-1/2}h^{-5/2}$, and (42) is satisfied with $c = \sigma a^{-5/2}$ if $d_n = h = an^{-1/5}$. From (32), it follows that Assumption A1 holds, with $w$ as in (43).

Define $g_n$ as in Assumption A2. Notice that

$$g_n(s) = h^{-2}\int l(u)\bar{m}_n(t_0 - hu)\,du$$
$$= h^{-2}\int l(u)m(t_0 - hu)\,du + r_n(s),$$

with $l(v) = k(v+s) - k(v) - sk'(v)$. Since

$$\int v^\lambda l(v)\,dv = \begin{cases} 0, & \lambda = 0,1, \\ s^2, & \lambda = 2, \end{cases}$$

it follows by a Taylor expansion of $m$ around $t_0$ that the first term converges to $As^2$, since $A = m''(t_0)/2$. The convergence is uniform with respect to $s$ over compact intervals, since the limit function $As^2$ is convex; see [25] and

Theorem 10.8 in [41]. For the second term, notice that

$$\sup_{|s|\leq c} |r_n(s)| \leq h^{-2} \sup_{|s|\leq c} \int |l(u)|\,du \sup_{|u-t_0|\leq(c+1)h} |\bar{m}_n(u) - m(u)|$$

$$= O(n^{-1}h^{-2}) = o(1)$$

since $nh^{-2} = (1+o(1))ch/\sigma_{\hat{n}} \to 0$, because $h \to 0$ and $\sigma_{\hat{n}} \to \infty$, which proves Assumption A2.

Assumptions A3 and A4 are satisfied by Proposition 1 and Lemmas B.1 and B.4 in Appendix B, and Assumption A5 holds by properties of Brownian motion. An application of Theorem 1 shows that

$$(44) \qquad d_n^{-2}(T_{[0,1]}(x_n)(t_0) - m(t_0)) \xrightarrow{\mathcal{L}} T(\tfrac{1}{2}m''(t_0)s^2 + \tilde{v}(s;t_0))(0)$$

as $n \to \infty$. Furthermore,

$$d_n^{-2}(x_n(t_0) - x_{b,n}(t_0))$$

$$(45) \qquad = d_n^{-2}v_n(t_0) = d_n^{-2}(nh)^{-1}\sigma_n \int k'(u)\bar{w}_n(t - uh)\,du$$

$$= d_n^{-3}n^{-1}\sigma_{nh} \int k'(u)\bar{w}_{nh}(h^{-1}t - u)\,du$$

$$\xrightarrow{\mathcal{L}} c\int k'(u)w(-u)\,du,$$

$$(46) \qquad d_n^{-2}(x_{b,n}(t_0) - m(t_0)) \to \tfrac{1}{2}m''(t_0)\int u^2 k(u)\,du,$$

as $n \to \infty$. Since the process $\bar{w}_n$ in (45) is the same as in the definition of $\tilde{v}_n$, one can make the rescaling in (44) and (45) simultanously to get joint convergence of (44) and (45); together with (46), this proves the theorem for the independent data case.

(ii) (Mixing case) The proof is similar to the proof of (i), replacing $\sigma$ by $\kappa$.

(iii) (Long range dependent data case) We want to choose $d_n = h$ so that (42) is satisfied with $c = |\eta_r|a$. Since $\sigma_{\hat{n}}^2 = \eta_r^2(nd_n)^{2-rd}l_1(nd_n)$, we get

$$(47) \qquad |\eta_r|a = d_n^{-3}n^{-1}|\eta_r|(nd_n)^{1-rd/2}l_1(nd_n)^{1/2}$$

$$\iff d_n = n^{-rd/(4+rd)}l_2(n;a),$$

where $l_2$ is another function slowly varying at infinity, implicitly defined in (47). We check the assumptions of Theorem 1 similarly as for (i) and (ii). □

In practice, it can be preferable to normalize the estimator, as in (40). It is an interesting problem to study the asymptotics for the normalized estimator $\tilde{m}$; we conjecture the same rate of convergence to hold and note that the integrated mean square error is smaller for the corrected estimator.

3.2.2. *The case $d_n/h \to \infty$ as $n \to \infty$.* This subcase is the least interesting for us, since the limit processes essentially are white noise processes; see [3].

3.2.3. *The case $d_n/h \to 0$ as $n \to \infty$.* In this case we can state limit distributions for $T_{[0,1]}(x_n)$ centered around $x_n$, the bias term, however, is of a larger order and thus the estimator has no useful statistical consequence; see [3] for details.

3.3. *Kernel estimation followed by isotonic regression.* Suppose $(t_i, y_i)$ are pairs of data satisfying the relation (31). Assuming that $m$ is increasing, there is an alternative to doing isotonic regression. Define first

$$\tilde{y}_n = y_i, \qquad t_i - \frac{1}{2n} < t \leq t_i + \frac{1}{2n}, \qquad i = 1, \ldots, n,$$

as the piecewise constant interpolation of $\{y_i\}$. Similarly, we define $\tilde{m}_n$ and $\tilde{\varepsilon}_n$ from $\{m(t_i)\}$ and $\{\varepsilon_i\}$. Compute the Gasser–Müller kernel estimate (see [19]),

$$m_n(t) = h^{-1} \int k\left(\frac{t-u}{h}\right) \tilde{y}_n(u)\, du,$$

of $m$ and then do isotonic regression on the data $(t, m_n(t))_{0 \leq t \leq 1}$. We do isotonic regression according to $\hat{m}(t) = T_{[0,1]}(x_n)'(t)$, where

(48) $$x_n(t) = \int_{-\infty}^{t} m_n(u)\, du = \int K\left(\frac{t-u}{h}\right) \tilde{y}_n(u)\, du$$

is the primitive of $m_n$, and $K(t) = \int_{-\infty}^{t} k(u)\, du$. This is considered in [34], where the limit distribution is given for i.i.d. data and for the particular choice of bandwidth $h = n^{-1/5}$. In [34] the reverse scheme is treated also, that is, isotonic regression followed by smoothing, for i.i.d. data.

The deterministic and stochastic parts of $x_n$ are defined according to

$$x_n(t) = x_{b,n}(t) + v_n(t),$$

with

$$x_{b,n}(t) = \int K\left(\frac{t-u}{h}\right) \tilde{m}_n(u)\, du,$$

$$v_n(t) = \int K\left(\frac{t-u}{h}\right) \tilde{\varepsilon}_n(u)\, du$$

$$= \frac{\sigma_n}{n} \int k(u) w_n(t - uh)\, du.$$

Notice that

$$x'_{b,n}(t) = h^{-1} \int k\left(\frac{t-u}{h}\right) \tilde{m}_n(u)\, du$$

is increasing since $\tilde{m}_n$ is, and thus $x_{b,n}$ is convex. Notice further that, for the bias of $x'_{b,n}$, we have

$$\begin{aligned}
x'_{b,n}(t_0) - m(t_0) &= h^{-1} \int k\left(\frac{t-u}{h}\right)(\tilde{m}_n(u) - m(t_0))\, du \\
&= \frac{1}{2} m''(t_0) h^2 \int u^2 k(u)\, du + o(h^2) + O(n^{-1})
\end{aligned} \tag{49}$$

as $n \to \infty$, if we assume that $m''(t_0)$ exists.

For the stochastic part, again we get different results depending on the asymptotic size of $d_n/h$.

3.3.1. *The case $d_n = h$.* The random part can, for $s > 0$, be rescaled as

$$\begin{aligned}
\tilde{v}_n(s;t) &= d_n^{-2} n^{-1} \sigma_{\hat{n}} \int (w_{\hat{n}}(h^{-1}t + s - u) - w_{\hat{n}}(h^{-1}t - u)) k(u)\, du \\
&\stackrel{\mathcal{L}}{=} d_n^{-2} n^{-1} \sigma_{\hat{n}} \int [w_{\hat{n}}(s-u) - w_{\hat{n}}(-u)] k(u)\, du,
\end{aligned}$$

with $\hat{n} = nd_n$, the right-hand side valid also for $s < 0$, and the last equality being exact only for $t = t_i$ and holding asymptotically otherwise. Assuming that

$$c_n = d_n^{-2} n^{-1} \sigma_{\hat{n}} \to c > 0, \tag{50}$$

the integrability of $k$, $w_n \stackrel{\mathcal{L}}{\to} w$ on $D(-\infty, \infty)$ and the continuous mapping theorem imply that

$$\tilde{v}_n(s;t) \stackrel{\mathcal{L}}{\to} \tilde{v}(s;t) = c \int (w(s-u) - w(-u)) k(u)\, du \tag{51}$$

on $D(-\infty, \infty)$ as $n \to \infty$.

Note that (49) implies that in the following two theorems $x'_{b,n}(t_0)$ can be replaced by $m(t_0)$.

THEOREM 5. *Assume $m$ is increasing, $m'(t_0) > 0$ and $t_0 \in (0,1)$. Assume that $\{\varepsilon_i\}$ are independent and identically distributed with $E(\varepsilon_i) = 0$ and $\mathrm{Var}(\varepsilon_i) = \sigma^2$. Define $x_n$ as in (48), with a nonnegative and compactly supported kernel $k$ having a bounded derivative $k'$, and with bandwidth $h$ specified below. Let $\hat{m}(t) = T_{[0,1]}(x_n)'(t)$. Suppose that one of the following conditions holds:*

(i) *$\{\varepsilon_i\}$ are independent and identically distributed with $E\varepsilon_i = 0$ and $\sigma^2 = \mathrm{Var}(\varepsilon_i) < \infty$ and we choose $h = an^{-1/3}$, where $a > 0$ is an arbitrary constant;*

(ii) *Assumption A8 or A9 holds, $\sigma_n^2 = \mathrm{Var}(\sum_{i=1}^n \varepsilon_i)$ and $\kappa^2$ is defined in (33) and we choose $h = an^{-1/3}$, where $a > 0$ is an arbitrary constant;*

(iii) $\varepsilon_i = g(\xi_i)$ is a long range dependent subordinated Gaussian sequence with parameters $d$ and $r$, and $\beta$ as in (36), and we choose $h = l_2(n;a)n^{-rd/(2+rd)}$, where $a > 0$ and $n \mapsto l_2(n;a)$ is a slowly varying function defined in the proof below.

Then, correspondingly, we obtain

$$d_n^{-1}(\hat{m}(t_0) - x'_{b,n}(t_0)) \xrightarrow{\mathcal{L}} (2m'(t_0))^{1/2} a \cdot \arg\min_{s \in \mathbb{R}} \left( s^2 + \tilde{v}\left(\frac{s}{\sqrt{m'(t_0)/2}}; t_0\right) \right),$$

as $n \to \infty$, where, respectively:

(i) $d_n = n^{-1/3}$, $w = B$ in (51), $c = a^{-3/2}\sigma$,
(ii) $d_n = n^{-1/3}$, $w = B$ in (51), $c = a^{-3/2}\kappa$,
(iii) $d_n = l_2(n;a)n^{-rd/(2+rd)}$, $w = B_{r,\beta}$ in (51), $c = |\eta_r|a$,

and with $c = a^{-3/2}\sigma, w = B$ in (51).

PROOF. (i) (Independent data case) Since $\sigma_{\hat{n}}^2 = \sigma^2\hat{n} = \sigma^2 n d_n$, we get $d_n^{-2}n^{-1}\sigma_{\hat{n}} = d_n^{-3/2}n^{-1/2}\sigma$. Putting $d_n = h = an^{-1/3}$, we thus get $c = a^{-3/2}\sigma$. Let us now verify the theorem from Corollary 1. From (32) and (51), we deduce Assumption A1. Notice that

$$x_{b,n}(t) = \int k(u)\bar{x}_{b,n}(t - uh) \, du,$$

with $\bar{x}_{b,n}(t) = \int_{1/2}^{t} \tilde{m}_n(u) \, du$ a piecewise linear approximation of the convex function $x_b(t) = \int_{1/2}^{t} m(u) \, du$. Thus, the rest of the proof of Assumption A2 is similar to Theorem 4(i), replacing $\bar{m}_n$ and $m$ in Theorem 4(i) with $\bar{x}_{b,n}$ and $x_{b,n}$ respectively. Clearly, $|\bar{x}_{b,n} - x_{b,n}| = O(n^{-1})$ uniformly on compact subsets of $(0,1)$, since the same is true for $|\tilde{m}_n - m|$. Furthermore, Assumptions A3 and A4 follow by Proposition 1 and Lemmas B.1 and B.4 in Appendix B.

Since $\tilde{v}$ has stationary increments, Assumption A7 holds and Assumption A5 is motivated as in previous results. Furthermore, since $k'$ exists, $\tilde{v}$ is differentiable and thus, Assumption A6 holds (cf. Note 1). Corollary 1, with $A = m'(t_0)/2$, now implies the theorem.

(ii) (Mixing data case) For mixing data, $\sigma_{\hat{n}} \sim \kappa\hat{n}^{1/2}$. The rest of the proof is similar to the independent data case.

(iii) (Long range dependent data) Choose $d_n = h$ so that (50) is satisfied with $c_n = |\eta_r|a$. Since the variance is $\sigma_{\hat{n}}^2 = (nd_n)^{2-rd}\eta_r^2 l_1(nd_n)$, with $l_1$ a slowly varying function, we get

$$|\eta_r|a = d_n^{-2}n^{-1}(nd_n)^{1-rd/2}|\eta_r|l_1(nd_n)^{1/2}$$
(52)
$$\iff d_n = n^{-rd/(2+rd)}l_2(n;a),$$

where $l_2$ is another function slowly varying at infinity, implicitly defined by (52). Thus, with $h = d_n$, we obtain the theorem from Corollary 1. Assumptions A1–A7 are checked as in the last two theorems. $\square$

3.3.2. *The case $d_n/h \to \infty$ as $n \to \infty$.* Rescale the random part as

$$\tilde{v}_n(s;t) \stackrel{\mathcal{L}}{=} d_n^{-2} n^{-1} \sigma_{\hat{n}} \int (w_{\hat{n}}(s - uh/d_n) - w_{\hat{n}}(-uh/d_n)) k(u) \, du,$$

with $\hat{n} = nd_n$, holding exactly for $t = t_i$ and asymptotically for all $t$. Assume that

$$(53) \qquad d_n^{-2} n^{-1} \sigma_{\hat{n}} \to c$$

as $n \to \infty$. Write $\tilde{v}_n(\cdot; t) \stackrel{\mathcal{L}}{=} \Lambda_n w_{\hat{n}}$, with $\Lambda_n$ an operator $D(-\infty, \infty) \to D(-\infty, \infty)$. Then $\Lambda_n w \to \Lambda w = cw$ as $n \to \infty$ whenever $w \in C(-\infty, \infty)$. Thus, with $w$ the weak limit of $\{w_n\}$, the extended continuous mapping theorem implies

$$(54) \qquad \tilde{v}_n(s; t_0) \stackrel{\mathcal{L}}{\to} \tilde{v}(s; t_0) = c(w(s) - w(0)) = cw(s).$$

But (53) and (54) are identical to the results for isotonic regression in Section 3.1 and, thus, the limit distributions must be the same.

For the bandwidths, (53) entails $h \ll n^{-1/3}$ in the independent and weakly dependent cases, and $h \ll l_2(n) n^{-rd/(2+rd)}$ in the long range dependent case, with $l_2(n)$ the same slowly varying function as in Theorem 3(iii).

3.3.3. *The case $d_n/h \to 0$ as $n \to \infty$.* We rescale as

$$\tilde{v}_n(s;t) = d_n^{-2} n^{-1} \sigma_{\hat{n}} \int (w_{\hat{n}}(sd_n/h - u) - w_{\hat{n}}(-u)) k(u) \, du$$

$$= d_n^{-1} (nh)^{-1} \sigma_{\hat{n}} \int w_{\hat{n}}(u) \frac{k(sd_n/h - u) - k(-u)}{d_n/h} \, du,$$

with $\hat{n} = nh$. Assume that $(d_n nh)^{-1} \sigma_{\hat{n}} \to c > 0$, and that $k$ is of bounded variation. Then, since $h/d_n (k(sd_n/h - u) - k(-u)) \to k'(-u)s$ for $s$ in compact sets, we can use the extended continuous mapping theorem to obtain

$$\tilde{v}_n(s;t) \stackrel{\mathcal{L}}{\to} \tilde{v}(s;t) = c \int w(u) k'(-u) \, du \cdot s = c\omega \cdot s,$$

where $w$ is the weak limit of $\{w_n\}$. Here $\omega \in N(0, \int k^2(u) \, du)$ for independent and weakly dependent data, $\omega$ is Gaussian for long range dependent data with rank $r = 1$ and non-Gaussian for $r > 1$. Note that

$$T(s^2 + c\omega s) = s^2 + c\omega s,$$
$$T(s^2 + c\omega s)'(0) = c\omega.$$

When the limit process is $T(s^2 + \omega s)(0) = 0$, this implies that we should study the rescaling and choice of normalizing constants more carefully, in order to get a nontrivial limit. In order to keep things simple, we skip this and give proofs only for the regression function.

We prove limit results for the properly normalized difference $T(x_n)'(t) - x_{b,n}(t)$ only. Note the relation (49) for the asymptotic bias. For independent data and with the choice of bandwidth $h = an^{-1/5}$, we have $(nh)^{1/2} \sim h^2$ and, thus, the asymptotic bias is of the same size as the variance and is by (49) equal to $\frac{1}{2}a^2 n^{-2/5} m''(t_0) \int u^2 k(u)\, du$. This is consistent with results in [34], where the independent data part of the following theorem was first proved for the special case $h = an^{-1/5}$.

THEOREM 6. *Assume $m$ is increasing, with $m'(t_0) > 0$ and $t_0 \in (0,1)$. Assume that $\{\varepsilon_i\}$ are independent and identically distributed with $E(\varepsilon_i) = 0$ and $\mathrm{Var}(\varepsilon_i) = \sigma^2$. Define $x_n$ as in (48), with a nonnegative and compactly supported kernel of bounded variation and with bandwidth specified below. Let $\hat{m}(t) = T_{[0,1]}(x_n)'(t)$. Suppose that one of the following conditions holds:*

(i) *$\{\varepsilon_i\}$ are independent and identically distributed with $E\varepsilon_i = 0$ and $\sigma^2 = \mathrm{Var}(\varepsilon_i) < \infty$ and we choose $h \gg an^{-1/3}$;*

(ii) *Assumption A8 or A9 holds, $\sigma_n^2 = \mathrm{Var}(\sum_{i=1}^n \varepsilon_i)$ and $\kappa^2$ is defined in (33) and we choose $h \gg an^{-1/3}$;*

(iii) *$\varepsilon_i = g(\xi_i)$ is a long range dependent subordinated Gaussian sequence with parameters $d$ and $r$, and $\beta$ as in (36), and we choose $h \gg l_2(n) n^{-rd/(2+rd)}$ and $n \mapsto l_2(n)$ is a slowly varying function defined in the proof below.*

*Then, correspondingly, we obtain*

$$d_n^{-1}(\hat{m}(t_0) - x_{b,n}(t_0)) \xrightarrow{\mathcal{L}} Z$$

*as $n \to \infty$, where, respectively:*

(i) *$d_n = (nh)^{-1/2}, Z = N(0, \sigma^2 \int k^2(u)\, du)$,*

(ii) *$d_n = (nh)^{-1/2}, Z = N(0, \kappa^2 \int k^2(u)\, du)$,*

(iii) *$d_n = l_1(nh)^{1/2}(nh)^{-rd/2}, Z = |\eta_r| \int k'(-u) B_{r,\beta}(u)\, du$, and with $l_1$ the slowly varying function defined in (35).*

PROOF. (i) (Independent data case) Since we have $\sigma_{\hat{n}}^2 = nh\sigma^2$, we see that $(d_n nh)^{-1}\sigma_{\hat{n}} = d_n^{-1}(nh)^{-1/2}\sigma$ converges to $c = \sigma$ if we choose $d_n = (nh)^{-1/2}$. Then the condition $d_n/h \to 0$ holds if $h \gg n^{-1/3}$. So for the stochastic part of the estimator, we obtain from Theorem 2

$$d_n^{-1}(\hat{m}(t_0) - x_{b,n}(t_0)) \xrightarrow{\mathcal{L}} \sigma\omega,$$

provided all regularity conditions are checked. Assumptions A5 and A6 are trivially satisfied, since $\tilde{v}(\cdot, t_0)$ is linear and $w$ is a continuous random variable with a symmetric distribution. To prove Assumption A2, write $x_{b,n}$ as in the proof of Theorem 5(i). Then we have

$$g_n(s) = h^{-2}\int l(u)\bar{x}_{b,n}(t_0 - uh)\, du$$

$$= h^{-2} \int l(u) x_b(t_0 - uh) \, du + r_n(s),$$

with

$$l(v) = \frac{h^2}{d_n^2} \left[ k\left(v + s\frac{d_n}{h}\right) - k(v) - s\frac{d_n}{h} k'(v) \right]$$

satisfying

$$\int v^\lambda l(v) \, dv = \begin{cases} 0, & \lambda = 0, 1, \\ s^2, & \lambda = 2. \end{cases}$$

Noting that $\sup_{|s| \le c} |r_n(s)| = O(n^{-1} h^{-2}) = o(n^{-1/3}) = o(1)$ as $n \to \infty$, the rest of the proof of Assumption A2 proceeds as in the proof of Theorem 4(i). Assumption A1 follows from (32) and the extended continuous mapping theorem, as noted above. Proposition 1 and Lemmas B.1 and C.1 imply Assumptions A3 and A4, which ends the proof.

(ii) (Mixing data case) Now we have $\sigma_{\hat{n}}^2 / \hat{n} \sim \kappa^2$, and the rest of the proof proceeds as for independent data.

(iii) (Long range dependent data case) Here

$$(d_n nh)^{-1} \sigma_{\hat{n}} = (d_n nh)^{-1} (nh)^{1-rd/2} |\eta_r| l_1(nh)^{1/2}$$
$$= d_n^{-1} (nh)^{-rd/2} |\eta_r| l_1(nh)^{1/2}.$$

With $c = |\eta_r|$, we obtain

$$|\eta_r| = d_n^{-1} (nh)^{-rd/2} |\eta_r| l_1(nh)^{1/2} \iff d_n = l_1(nh)^{1/2} (nh)^{-rd/2}.$$

Thus,

$$l_1(nh)^{-1/2} (nh)^{rd/2} (\hat{m}(t_0) - x_{b,n}(t_0)) \xrightarrow{\mathcal{L}} \omega |\eta_r|,$$

where $\omega = \int k'(-u) B_{r,\beta}(u) \, du$. The condition $d_n/h \to 0$ is satisfied if we let $h \gg l_2(n) n^{-rd/(2+rd)}$, where $l_2$ is any of the slowly varying functions $l_2(\cdot; a)$ defined in the proof of Theorem 5(iii). Assumptions A1–A6 are checked as in parts (i) and (ii). □

Note that since $T(s^2 + cws)(0) = (s^2 + cws)|_{s=0} = 0$, it follows from Theorems 6 and 1 that

$$\int_0^{t_0} \hat{m}(s) \, ds - x_n(t_0)$$
$$= o_P(d_n^2) = \begin{cases} o_P((nh)^{-1}) & \text{(independent, weakly dependent data)}, \\ o_P(l_1(nh)(nh)^{-rd}) & \text{(long range dependent data)}. \end{cases}$$

**4. Density and distribution function estimation.** The empirical distribution function is of a fundamental importance in distribution and density function estimation. To define this, assume $\{t_i\}_{i=-\infty}^{\infty}$ is a stationary sequence of random variables with a marginal distribution function $F$. Then the empirical distribution function is

$$F_n(t) = \frac{1}{n} \sum_{i=1}^{n} \mathbb{1}_{\{t_i \leq t\}}.$$

Note that $F_n$ is right continuous with left-hand limits and, thus, $F_n$ lies in the space $D(-\infty, \infty)$.

Note also that

$$F_n(t) = F(t) + F_n^0(t),$$

where

$$F_n^0(t) = \frac{1}{n} \sum_{i=1}^{n} (\mathbb{1}_{\{t_i \leq t\}} - F(t))$$

is the centered empirical process. Consider a sequence $\delta_n$ such that $\delta_n \downarrow 0, n\delta_n \uparrow \infty$ as $n \to \infty$. Define the centered empirical process locally around $t_0$ on a scale $\delta_n$ as

$$
\begin{aligned}
w_{n,\delta_n}(s; t_0) &= \sigma_{n,\delta_n}^{-1} n(F_n^0(t_0 + s\delta_n) - F_n^0(t_0)) \\
&= \sigma_{n,\delta_n}^{-1} \sum_{i=1}^{n} (\mathbb{1}_{\{t_i \leq t_0 + s\delta_n\}} - \mathbb{1}_{\{t_i \leq t_0\}} - F(t_0 + s\delta_n) + F(t_0)),
\end{aligned}
$$
(55)

where

$$
\begin{aligned}
\sigma_{n,\delta_n}^2 &= \mathrm{Var}(n(F_n^0(t_0 + \delta_n) - F_n^0(t_0))) \\
&= \mathrm{Var}\left(\sum_{i=1}^{n} (\mathbb{1}_{\{t_0 < t_i \leq t_0 + \delta_n\}} - F(t_0 + \delta_n) + F(t_0))\right).
\end{aligned}
$$

We will prove weak convergence $w_{n,\delta_n} \xrightarrow{\mathcal{L}} w$, on $D(-\infty, \infty)$, as $n \to \infty$, for independent, weakly dependent and subordinated Gaussian long range dependent data $\{t_i\}$.

THEOREM 7. *Assume $\{t_i\}$ are independent, $f(t_0) = F'(t_0)$ exists and $\delta_n \downarrow 0$, $n\delta_n \uparrow \infty$ as $n \to \infty$. Then*

$$\frac{\sigma_{n,\delta_n}^2}{n\delta_n f(t_0)} \to 1$$

*and*

$$w_{n,\delta_n}(s; t_0) \xrightarrow{\mathcal{L}} B(s)$$

*on $D(-\infty, \infty)$, as $n \to \infty$, where $B$ is a standard Brownian motion.*

The proof of this theorem is a standard application of the Cramér–Wold device and tightness; see the technical report of Anevski and Hössjer [3].

For weakly dependent data, we will use mixing conditions. Define the $\sigma$-algebras

$$\mathcal{F}_k = \sigma\{t_i : i \leq k\},$$
$$\bar{\mathcal{F}}_k = \sigma\{t_i : i \geq k\}.$$

Then the sequence $\{t_i\}$ is said to be $\phi$-mixing if Definition 1 in Section 4 is applicable, with $t_i$ in place of $\varepsilon_i$.

THEOREM 8. *Assume the stationary sequence $\{t_i\}$ is $\phi$-mixing with $\sum_{i=1}^{\infty} i\phi^{1/2}(i) < \infty$, that $\delta_n \to 0, n\delta_n \to \infty$ as $n \to \infty$, $f(t_0) = F'(t_0)$ exists, as well as the joint density $f_k(s_1, s_2)$ of $(t_1, t_{1+k})$ on $[t_0 - \delta, t_0 + \delta]^2$ for some $\delta > 0$, and $k \geq 1$. Assume also that we have the bound*

$$\sum_{k=1}^{\infty} M_k < \infty,$$

*with $M_k = \sup_{t_0 - \delta \leq s_1, s_2 \leq t_0 + \delta} |f_k(s_1, s_2) - f(s_1)f(s_2)|$. Then*

$$\text{(56)} \qquad \frac{\sigma_{n,\delta_n}^2}{n\delta_n f(t_0)} \to 1$$

*and*

$$w_{n,\delta_n}(s; t_0) \xrightarrow{\mathcal{L}} B(s)$$

*on $D(-\infty, \infty)$, as $n \to \infty$.*

For a proof of this theorem, see the technical report of Anevski and Hössjer [3].

In the long range dependent case, as in the partial sum process case, we make an expansion in Hermite polynomials of the terms in the sum defining the empirical distribution function at $t \in \mathbb{R}$. In this case, however, the terms depend on $t$, which makes the analysis somewhat different.

Thus, assume $\{\xi_i\}_{i \geq 1}$ is a stationary Gaussian process with mean zero and covariance function $\text{Cov}(k) = E(\xi_i \xi_{i+k})$ such that $\text{Cov}(0) = 1$ and $\text{Cov}(k) = k^{-d} l_0(k)$, where $l_0$ is a function slowly varying at infinity and $0 < d < 1$ is fixed. Let $g : \mathbb{R} \mapsto \mathbb{R}$ be a measurable function and $t_i = g(\xi_i)$. For a fixed $t$, expand the function $\mathbb{1}_{\{t_1 \leq t\}} - F(t)$ in Hermite polynomials

$$\mathbb{1}_{\{t_1 \leq t\}} - F(t) = \sum_{k=r(t)}^{\infty} \frac{1}{k!} \eta_k(t) h_k(\xi_1).$$

Here $h_k$ is the Hermite polynomial of order $k$, and

$$\eta_k(t) = E[(\mathbb{1}_{\{t_1 \leq t\}} - F(t))h_k(\xi_1)] = \int (\mathbb{1}_{\{g(u) \leq t\}} - F(t))h_k(u)\phi(u)\,du$$

are the $L^2(\phi)$-projections on $h_k$, and $r(t)$ is the first nonzero coefficient in the expansion. Now let $t$ vary and define the Hermite rank of the functions $\{\mathbb{1}_{\{g(\cdot) \leq t\}} - F(t) : t \in \mathbb{R}\}$ as $r = \inf_t r(t)$. Assume that $0 < d < 1/r$. With a slight abuse of notation, we say that the sequence $\{t_i\}$ is long range dependent subordinated Gaussian with parameters $d$ and $r$.

This implies that the sequence $\{\mathbb{1}_{\{t_i \leq t\}} - F(t)\}_{i \geq 1}$ exhibits long range dependence, and $\sigma_n^2 = \text{Var}(\sum_{i=1}^n h_r(\xi_i))$ is asymptotically proportional to $n^{2-rd}l_1(n)$, with $l_1$ defined in (35). From Theorem 1.1 in [17], under the above assumptions it follows that

$$\sigma_n^{-1} n F_n^0(t) \xrightarrow{\mathcal{L}} \frac{\eta_r(t)}{r!} z_r$$

on $D[-\infty, \infty]$ equipped with the supnorm-metric. The random variable $z_r$ is the evaluation $z_r = z_r(1)$ of the process defined in Section 4, with $\beta$ as in (36). Note that $z_r$ is Gaussian for $r = 1$ and non-Gaussian for $r \geq 2$. Note also that the space here is the compact $D[-\infty, \infty]$ and the metric is the supnorm-metric over the whole extended real line.

THEOREM 9. *Assume $\{t_i\}$ is a long range dependent subordinated Gaussian sequence with parameters $d$ and $r$ and $0 < d < 1/r$. Define*

$$\kappa_1 = \min(d, 1 - rd)/2,$$
$$\kappa_2 = \min(2d, 1 - rd)/2.$$

*Assume that $\delta_n \to 0$ as $n \to \infty$ and, for some $\varepsilon > 0$,*

$$\delta_n \gg n^{-\kappa_1 + \varepsilon} \qquad \text{if } d \geq 1/(1+r),$$
$$\delta_n \gg n^{-\kappa_2 + \varepsilon} \qquad \text{if } 0 < d < 1/(1+r).$$

*Then if $\eta_r$ and $F$ are differentiable at $t_0$ with $\eta_r'(t_0) \neq 0$,*

(57) $$\sigma_{n,\delta_n}^2 \sim \sigma_n^2 \left(\frac{\eta_r'(t_0)}{r!}\right)^2 \delta_n^2$$

*and*

(58) $$w_{n,\delta_n}(s;t_0) \xrightarrow{\mathcal{L}} s \cdot sgn(\eta_r'(t_0)) z_r$$

*as $n \to \infty$, on $D(-\infty, \infty)$.*

PROOF. Write

$$nF_n^0(t) = \eta_r(t) \sum_{i=1}^{n} \frac{h_r(\xi_i)}{r!} + \sigma_n S_n(t),$$

with $S_n$ containing the higher-order terms in the Hermite expansion,

$$S_n(t) = \sigma_n^{-1} \sum_{k=r+1}^{\infty} \frac{1}{k!}\eta_k(t) \sum_{i=1}^{n} h_k(\xi_i).$$

Then

$$n(F_n^0(t_0 + \delta_n) - F_n^0(t_0)) = \delta_n(\eta_r'(t_0) + o(1)) \sum_{i=1}^{n} \frac{h_r(\xi_i)}{r!}$$
$$+ \sigma_n(S_n(t_0 + \delta_n) - S_n(t_0))$$

as $n \to \infty$. Thus, to prove (57) it suffices to show that

$$\text{Var}(S_n(t_0 + \delta_n) - S_n(t_0)) = o(\delta_n^2)$$

as $n \to \infty$. With $\tilde{\sigma}_n^2 = \text{Var}(\sum_{i=1}^{n} h_{r+1}(\xi_i))$ we get, for large enough $n$,

$$\text{Var}(S_n(t_0 + \delta_n) - S_n(t_0)) \leq 2(\text{Var}(S_n(t_0 + \delta_n)) + \text{Var}(S_n(t_0)))$$
$$\sim 4\eta_{r+1}(t_0)^2(\tilde{\sigma}_n/\sigma_n)^2 \leq n^{-2\kappa_1 + \varepsilon}$$

proving (57) if $d \geq 1/(1+r)$, since then, by assumption, $\delta_n \gg n^{-\kappa_1 + \varepsilon}$.

If instead $d < 1/(1+r)$, we define $\tilde{S}_n$ by

$$\sigma_n S_n(t) = \frac{\eta_{r+1}(t)}{(r+1)!} \sum_{i=1}^{n} h_{r+1}(\xi_i) + \tilde{\sigma}_n \tilde{S}_n(t)$$

and, thus,

$$|S_n(t_0 + s\delta_n) - S_n(t_0)|$$
(59)
$$\leq \frac{\tilde{\sigma}_n}{\sigma_n} \left( \frac{\sum_{i=1}^{n} h_{r+1}(\xi_i)}{\tilde{\sigma}_n(r+1)!} |\eta_{r+1}(t_0 + s\delta_n) \right.$$
$$\left. - \eta_{r+1}(t_0)| + |\tilde{S}_n(t_0 + s\delta_n) - \tilde{S}_n(t_0)| \right).$$

Using the relations (cf. page 997 of [12]),

(60) $\quad (\eta_{r+1}(t_0 + \delta_n) - \eta_{r+1}(t_0))^2 \leq (r+1)!(F(t_0 + \delta_n) - F(t_0)),$

(61) $\quad \tilde{\sigma}_n^2 \leq \sigma_n^2 n^{-d+\varepsilon},$

$$\text{Var}(\tilde{S}_n(t)) \leq n^{-2\kappa_2 + d + \varepsilon},$$

which hold for large enough $n$, from (59) we get

$$\mathrm{Var}(S_n(t_0 + \delta_n) - S_n(t_0))$$
$$\leq \left(\frac{\tilde{\sigma}_n}{\sigma_n}\right)^2 O[\delta_n^2 + 4(\mathrm{Var}(\tilde{S}_n(t_0 + \delta_n)) + \mathrm{Var}(\tilde{S}_n(t_0)))]$$
$$= O(n^{-d+\varepsilon}\delta_n^2 + n^{-2\kappa_2+\varepsilon}) = o(\delta_n^2)$$

as $n \to \infty$, since, by assumption, $\delta_n \gg n^{-\kappa_2+\varepsilon}$.

To prove (58), notice that

$$w_{n,\delta_n}(s; t_0) = \sigma_{n,\delta_n}^{-1} n(F_n^0(t_0 + s\delta_n) - F_n^0(t_0))$$
$$= (1 + o(1))\frac{\eta_r(t_0 + s\delta_n) - \eta_r(t_0)}{\delta_n|\eta_r'(t_0)|}\sigma_n^{-1}\sum_{i=1}^n h_r(\xi_i)$$
$$+ C_n \delta_n^{-1}(S_n(t_0 + s\delta_n) - S_n(t_0)),$$

where $C_n \to r!/|\eta_r'(t_0)|$, as $n \to \infty$. Since $\sigma_n^{-1}\sum_{i=1}^n h_r(\xi_i) \overset{\mathcal{L}}{\to} z_r$, and

$$\frac{\eta_r(t_0 + s\delta_n) - \eta_r(t_0)}{\delta_n|\eta_r'(t_0)|} \to s \cdot \mathrm{sgn}(\eta_r'(t_0))$$

uniformly on compacts as $n \to \infty$, (58) will follow if we establish

(62) $$\sup_{s \in \mathbb{R}} \delta_n^{-1}|S_n(t_0 + s\delta_n) - S_n(t_0)| \overset{P}{\to} 0$$

as $n \to \infty$.

If $d \geq 1/(1+r)$, then from formula (2.2) in [12] we obtain

$$\delta_n^{-1}|S_n(t_0 + s\delta_n) - S_n(t_0)| \leq 2\delta_n^{-1}\sup_{t \in \mathbb{R}}|S_n(t)|$$
$$= O_P(\delta_n^{-1}n^{-\kappa_1+\varepsilon}),$$

which proves (62), since $\delta_n \gg n^{-\kappa_1+\varepsilon}$. If $d < 1/(1+r)$, then from (59), (60), (61) and formula (2.3) in [12], we have

$$\delta_n^{-1}|S_n(t_0 + s\delta_n) - S_n(t_0)| \sim \delta_n^{-1} n^{-d/2+\varepsilon}(\delta_n^{1/2} + O_P(n^{-\kappa_2+d/2+\varepsilon}))$$

as $\delta_n \to 0$. Since $\delta_n \gg n^{-\kappa_2+\varepsilon}$, implying also $\delta_n \gg n^{-d+\varepsilon}$, (62) holds and thus, (58) is proved. $\square$

4.1. *Estimating an increasing density function.* Suppose we have observations from an unknown density $f$ lying in the class $\mathcal{F} = \{f : (-\infty, 0] \mapsto [0, \infty), f \geq 0, \int f(u)\, du = 1, f \text{ increasing}\}$, and assume we want to estimate $f$ at a fixed point $t$. In the case of independent data, we can easily write down the likelihood and try to maximize this over the class $\mathcal{F}$. The solution is the nonparametric maximum likelihood estimate, and it is known to be given by $T_{(-\infty,0]}(F_n)'$ (see [20]), where $F_n$ is the empirical distribution function. In the case of independent data, also the limit distribution of $T_{(-\infty,0]}(F_n)$ and $T_{(-\infty,0]}(F_n)'$ are known; see [21, 22, 39] and [46]. We will put these results into a more general framework.

The algorithm $T_{(-\infty,0]}(F_n)'$ produces an increasing density also in the case of dependent data, with marginal $f$, while of course the likelihood is more difficult to work with. Thus, $T_{(-\infty,0]}(F_n)'$ is an ad hoc estimator of an increasing density in the case of dependent data which lies in $\mathcal{F}$, and for which we will derive the limit distribution.

Let

$$x_n(t) = F_n(t),$$
$$x_{b,n}(t) = F(t),$$
$$v_n(t) = F_n^0(t).$$

Under various dependence assumptions, we have

$$\tilde{v}_n(s; t_0) = c_n w_{n,d_n}(s; t_0) \xrightarrow{\mathcal{L}} w(s) =: v(s; t_0)$$

on $D(-\infty, \infty)$, as $n \to \infty$, where $\{d_n\}$ is chosen so that

$$c_n = d_n^{-2} n^{-1} \sigma_{n,d_n} \to 1.$$

THEOREM 10 (Independent and mixing data). *Assume $\{t_i\}_{i \geq 1}$ is a stationary sequence with an an increasing marginal density function $f$ such that $f'(t_0) > 0$ and $t_0 < 0$. Let $F_n(t)$ be the empirical distribution function and $\hat{f}_n(t) = T_{(-\infty,0]}(F_n)'(t)$. Suppose that one of the following conditions holds:*

  (i) *$\{t_i\}_{i \geq 1}$ is an i.i.d. sequence;*
  (ii) *$\{t_i\}_{i \geq 1}$ satisfies the assumptions of Theorem 8.*

*Then we obtain*

$$n^{1/3} c_1(t_0)(\hat{f}_n(t_0) - f(t_0)) \xrightarrow{\mathcal{L}} \arg\min_{s \in \mathbb{R}}(s^2 + B(s)),$$

$$n^{2/3} c_2(t_0)\left(\int_0^{t_0} \hat{f}_n(s)\,ds - F(t_0)\right) \xrightarrow{\mathcal{L}} T(s^2 + B(s))(0)$$

*as $n \to \infty$, with*

$$c_1(t_0) = f(t_0)^{-1/3}(\tfrac{1}{2})^{2/3} f'(t_0)^{-1/3},$$
$$c_2(t_0) = f(t_0)^{-2/3}(\tfrac{1}{2} f'(t_0))^{1/3},$$

*and $B$ a standard two-sided Brownian motion.*

PROOF. (i) (Independent data case) To determine the constants, we use Theorem 7,

$$\sigma_{n,d_n}^2 \sim nd_n f(t_0)$$

$$\iff c_n = d_n^{-2} n^{-1} \sqrt{nd_n f(t_0)} = d_n^{-3/2} n^{-1/2} f(t_0)^{1/2} \sim 1$$

$$\iff d_n \sim f(t_0)^{1/3} n^{-1/3}.$$

Since $x_{b,n} = F$ is convex and $x_b''(t_0) = f'(t_0)$, Assumption A2 is satisfied with $A = \frac{1}{2} f'(t_0)$. From Theorem 7, it follows that Assumption A1 is satisfied, and Proposition 1 and Lemma C.1 imply Assumptions A3 and A4. Assumptions A5, A6 and A7 hold by properties of Brownian motion (cf. [42] for an LIL for Brownian motion implying Assumption A6 by Note 1), so that Corollary 1 implies

$$d_n^{-1}(\hat{f}_n(t_0) - f(t_0)) \xrightarrow{\mathcal{L}} 2\sqrt{A} \arg\min_{s \in \mathbb{R}} \left( s^2 + B\left(\frac{s}{\sqrt{A}}\right) \right)$$

$$\stackrel{\mathcal{L}}{=} 2A^{1/3} \arg\min_{s \in \mathbb{R}} (s^2 + B(s))$$

as $n \to \infty$ and, thus, $c_1(t_0)$ has the form stated in the theorem. Then Theorem 1 implies

$$d_n^{-2} \left( \int_0^{t_0} \hat{f}_n(s) \, ds - F(t_0) \right) \xrightarrow{\mathcal{L}} T(As^2 + B(s))(0)$$

$$\stackrel{\mathcal{L}}{=} A^{-1/3} T(s^2 + B(s))(0)$$

as $n \to \infty$ and, thus, $c_2(t_0)$ has the form stated in the theorem.

(ii) (Mixing data case) The proof is completely analoguous to the i.i.d. case and uses Theorem 8 instead of Theorem 7. □

THEOREM 11 (Long range dependent data). *Assume $\{t_i\}_{i \geq 1}$ is a long range dependent subordinated Gaussian sequence with parameters $r = 1$ and $0 < d < 1/2$, and $\beta$ as in (36). Let $f$ be the marginal density function of $\{t_i\}$, and assume $f$ is increasing with $f'(t_0) > 0$ and $t_0 < 0$. Then with $F_n$ the empirical distribution function and $\hat{f}_n(t) = T_{(-\infty,0]}(F_n)'(t)$,*

$$|\eta_1'(t_0)|^{-1} l_1(n)^{-1/2} n^{d/2} (\hat{f}_n(t_0) - f(t_0)) \xrightarrow{\mathcal{L}} N(0,1)$$

*as $n \to \infty$, with $l_1$ a function slowly varying at infinity, defined as in (35).*

PROOF. We have

$$c_n = d_n^{-2} n^{-1} \sigma_{n,d_n} \sim d_n^{-1} n^{-1} \sigma_n |\eta_r'(t_0)|/r! \sim 1$$

$$\iff d_n \sim |\eta_r'(t_0)| l_1(n)^{1/2} n^{-rd/2}/r!,$$

where $\sigma_n$ and $l_1$ are defined before Theorem 9. Note that the assumptions in Theorem 9, with $\delta_n = d_n$, are only satisfied for $r = 1, 0 < d < 1/2$, for which case we have $w_{n,d_n} \xrightarrow{\mathcal{L}} s \cdot z_1 =: w(s)$ as $n \to \infty$, with $z_1$ a standard Gaussian random variable. Theorem 2 implies

$$d_n^{-1} \hat{f}_n((t_0) - f(t_0)) \xrightarrow{\mathcal{L}} T(s^2 + z_1 s)'(0) = z_1$$

as $n \to \infty$, implying the theorem. Assumption A1 follows from Theorem 9, Assumption A2 is established as in Theorem 10(i), Assumptions A3 and A4 follow from Proposition 1 and Lemma C.1 and Assumptions A5 and A6 are trivially satisfied; see Note 1. $\square$

4.2. *Estimating a convex density function.* Suppose $f:[0,\infty) \mapsto [0,\infty)$ is a convex density function, and we want to find an estimator of $f$ at a fixed point $t_0 > 0$. For independent data, it is possible to define the nonparametric maximum likelihood estimate. The algorithm for calculating this is quite complicated though; see [30]; see also [23, 24] for the limit distribution. We present the following alternative estimator: Let

$$x_n(t) = n^{-1} h^{-1} \sum_{i=1}^{n} k\left(\frac{t - t_i}{h}\right)$$

be the kernel estimator of the density $f$, with $k$ a density function supported on $[-1, 1]$, and $h > 0$ the bandwidth. Define the (nonnormalized) density estimate $\hat{f}_n(t) = T(x_n)(t)$, and note that $\hat{f}_n$ is convex and positive, but does not integrate to one [note that the estimator $T_{[0,\infty)}(x_n)/I_n$ is a convex density function, where $I_n = \int T_{[0,\infty)}(x_n)(s)\,ds$]. We will state the limit distributions for $\hat{f}_n$ in the weakly dependent cases; see Section 1 for an interpretation. In the long range dependent case the limit process is pathological $[\tilde{v}(s;t_0) = 0]$, so that the rate of convergence is faster than indicated by our approach. Since a further study of this case is not straightforward, we refrain from more work on this; see the remark after Theorem 12.

We can write

$$x_n(t) = h^{-1} \int k'(u) F_n(t - hu)\,du,$$

$$x_{b,n}(t) = h^{-1} \int k'(u) F(t - hu)\,du,$$

$$v_n(t) = h^{-1} \int k'(u) F_n^0(t - hu)\,du.$$

4.2.1. *The case $d_n = h$.* The rescaled process is

$$\tilde{v}_n(s; t_0) = c_n \int k'(u)(w_{n,d_n}(s - u; t_0) - w_{n,d_n}(-u; t_0))\,du,$$

with $c_n = d_n^{-2}(nh)^{-1}\sigma_{n,d_n}$. Choosing $d_n$ so that $c_n \to c$ as $n \to \infty$, for some constant $c$, we obtain

$$\tilde{v}_n(s;t_0) \xrightarrow{\mathcal{L}} c \int k'(u)(w(s-u) - w(-u))\, du =: \tilde{v}(s;t_0) \tag{63}$$

on $D(-\infty, \infty)$ as $n \to \infty$, using the continuous mapping theorem as in Section 4.2.1, and with $w$ the weak limit of $\{w_n\}$, assuming that $k'$ is bounded and since $k$ has compact support.

Recall the definition of the (nonnormalized) density estimate $\hat{f}_n(t) = T(x_n)(t)$. The rate $n^{-2/5}$ for the estimator $\hat{f}_n$ in Theorem 12 is the same as in the limit distribution of the NPMLE; see [23, 24]. This is also the optimal rate for estimating a convex density from independent observations; see [2].

THEOREM 12 (Independent and mixing data). *Let $\{t_i\}_{i \geq 1}$ be a stationary sequence with a convex marginal density function $f$ such that $f''(t_0) > 0$ and $t_0 > 0$. Let $x_n(t)$ be the kernel density function above with $k$ a compactly supported density such that $k'$ is bounded, $h = an^{-1/5}$ and $a > 0$ an arbitrary constant. Suppose that one of the following conditions holds:*

(i) *$\{t_i\}_{i \geq 1}$ is an i.i.d. sequence;*
(ii) *$\{t_i\}_{i \geq 1}$ satisfies the assumptions of Theorem 9.*

*Then we obtain*

$$n^{2/5}(\hat{f}_n(t_0) - f(t_0))$$
$$\xrightarrow{\mathcal{L}} a^2 T(\tfrac{1}{2}f''(t_0)s^2 + \tilde{v}(s;t_0))(0)$$
$$+ \tfrac{1}{2}a^2 \int u^2 k(u)\, du\, f''(t_0) + ca^2 \int k'(u)w(-u)\, du$$

*as $n \to \infty$, with $c = a^{-5/2}f(t_0)^{1/2}$, $\tilde{v}(s;t)$ as in (63) and $w$ a standard two-sided Brownian motion.*

PROOF. (i) (Independent data case) We have $\sigma_{n,d_n}^2 \sim nd_n f(t_0)$, so that

$$d_n^{-2}(nh)^{-1}\sigma_{n,d_n} \sim d_n^{-5/2} n^{-1/2} f(t_0)^{1/2}.$$

If $d_n = an^{-1/5}$, we get $c = a^{-5/2}f(t_0)^{1/2}$ and

$$n^{2/5}(\hat{f}_n(t_0) - f(t_0)) \xrightarrow{\mathcal{L}} a^2 T(\tfrac{1}{2}f''(t_0)s^2 + \tilde{v}(s;t_0))(0) \tag{64}$$

follows from Theorem 1, provided the conditions in Theorem 1 hold. Notice that

$$x_{b,n}(t) = h^{-1} \int k\left(\frac{t-u}{h}\right) f(u)\, du$$

is convex, which establishes Assumption A2, with $A = \frac{1}{2}f''(t_0)$, similar to the argument in Section 4.2.1. Assumptions A1, A3 and A4 and A5 are also verified analogously to the argument in Section 4.2.1. Finally,

$$
(65) \quad \begin{aligned} d_n^{-2}(x_n(t_0) - x_{b,n}(t_0)) &= d_n^{-3} n^{-1} \sigma_{nd_n} \int k'(u) w_{nd_n}(-u; t_0) \, du \\ &\stackrel{\mathcal{L}}{\to} c \int k'(u) w(-u) \, du \end{aligned}
$$

$$
(66) \quad d_n^{-2}(x_{b,n}(t_0) - f(t_0)) \to \tfrac{1}{2} \int u^2 k(u) \, du f''(t_0).
$$

The joint convergence in (64) and (65) (cf. the proof of Theorem 4), together with (66), shows the statement of the theorem for the independent data case.

(ii) (Mixing data case) Similar to the proof of case (i). □

For long range dependent data, as in Section 4.1, we are restricted to the case $r = 1, 0 < d < 1/2$. But now

$$
\begin{aligned} \tilde{v}(s; t_0) &= c \int k'(u)((s-u)z_r - (-u)z_r) \, du \\ &= csz_r \int k'(u) \, du = 0, \end{aligned}
$$

where the first equality holds since $w(s) \stackrel{\mathcal{L}}{=} s \cdot z_1$. This indicates that the rate of convergence is faster than obtained solving $c_n = d_n^{-3} n^{-1} \sigma_{n,d_n} = c$. We refrain from further work for this case.

4.2.2. *The cases $d_n/h \to \infty$ and $d_n/h \to \infty$ as $n \to \infty$.* We refer the interested reader to the technical report of Anevski and Hössjer [3].

**5. Self similarity and rates of convergence.** In many of the examples treated in Sections 3 and 4, the stochastic part $v_n$ of $x_n$ is asymptotically self similar in the following sense: There exists a sequence $a_n \downarrow 0$ such that $a_n^{-1} v_n$ converges in distribution on a scale $b_n$ around $t_0$. More precisely, we assume the existence of a limit process $\bar{v}(\cdot; t_0)$ such that

$$
(67) \quad a_n^{-1} v_n(t_0 + sb_n) \stackrel{\mathcal{L}}{\to} \bar{v}(s; t_0).
$$

For local estimators, $b_n \downarrow 0$ and the convergence in (67) takes place in $D(-\infty, \infty)$. For global estimators, we have $b_n \equiv 1$ and then the convergence takes place in $D(J - t_0)$. Further, assume that $\bar{v}(\cdot; t_0)$ is locally self similar, in the sense that, for some $\beta > 0$ and some process $\tilde{v}(\cdot; t_0)$,

$$
(68) \quad \delta^{-\beta}(\bar{v}(\delta s; t_0) - \bar{v}(0; t_0)) \stackrel{\mathcal{L}}{\to} \tilde{v}(s; t_0)
$$

TABLE 1
*Convergence rates $d_n$ for various choices of $a_n$, $b_n$ and $\beta$*

| Theorem | $a_n$ | $b_n$ | $\beta$ | $d_n$ |
|---|---|---|---|---|
| 3(i), (ii), 10(i), (ii) | $n^{-1/2}$ | 1 | $1/2$ | $n^{-1/3}$ |
| 3(iii) | $n^{-rd/2}$ | 1 | $1 - rd/2$ | $n^{-rd/(2+rd)}$ |
| 6(i), (ii) | $n^{-1/2}h^{1/2}$ | $h$ | 1 | $(nh)^{-1/2}$ |
| 6(iii) | $n^{-rd/2}h^{1-rd/2}$ | $h$ | 1 | $(nh)^{-rd/2}$ |
| 11 | $n^{-d/2}$ | 1 | 1 | $n^{-d/2}$ |

on $D(-\infty, \infty)$ as $\delta \to 0$. Suppose that $d_n \ll b_n$ and put $\delta_n = d_n/b_n$. If we can interchange limits between (67) and (68), we obtain

$$
\begin{aligned}
(69) \quad d_n^{-p}(v_n(t_0 + sd_n) - v_n(t_0)) &\stackrel{\mathcal{L}}{\approx} d_n^{-p} a_n(\bar{v}(s\delta_n; t_0) - \bar{v}(0; t_0)) \\
&\stackrel{\mathcal{L}}{\approx} d_n^{-p} a_n \delta_n^{\beta} \tilde{v}(s; t_0).
\end{aligned}
$$

Thus, Assumption A1 requires [up to a factor $1 + o(1)$]

$$
(70) \qquad d_n^{-p} a_n \delta_n^{\beta} = 1 \quad \Longleftrightarrow \quad d_n = (a_n b_n^{-\beta})^{1/(p-\beta)},
$$

which is then a general formula for choosing $d_n$. As a consequence of this

$$
T(x_n)(t_0) - x_n(t_0) = O_P(d_n^p) = O_P((a_n b_n^{-\beta})^{p/(p-\beta)}),
$$
$$
T(x_n)'(t_0) - x'_{b,n}(t_0) = O_P(d_n^{p-1}) = O_P((a_n b_n^{-\beta})^{(p-1)/(p-\beta)}).
$$

For instance, in Theorem 3(i) we have, by Donsker's theorem, $a_n = n^{-1/2}$, $b_n = 1$ and $\bar{v}(s; t_0) = \sigma B(s + t_0), s \in J - t_0$, with $B$ a standard Brownian motion. Since the Brownian motion has stationary increments and is self similar with $\beta = 1/2$, we can put $\tilde{v}(s; t_0) \stackrel{\mathcal{L}}{=} \sigma B(s)$ on $D(-\infty, \infty)$ so that

$$
d_n = a_n^{1/(p-\beta)} = n^{-1/2(1/(p-1/2))} = n^{-1/(2p-1)},
$$

that is, $d_n = n^{-1/3}$ when $p = 2$.

Table 1 lists values of $a_n, b_n, \beta$ and $d_n$ for all examples with $d_n \ll b_n$ in Sections 3 and 4 when $p = 2$. For the long range dependent examples, we have simplified and put $l_1(n) = 1$ for the slowly varying function $l_1$ in (35). For general $l_1$, formula (67) is not valid. We have also ignored constants (not depending on $n$) of $d_n$.

For Theorem 6, we have replaced (67) with the more general requirement

$$
(71) \qquad a_n^{-1}(v_n(t_0 + sb_n) - v_n(t_0)) \stackrel{\mathcal{L}}{\to} \bar{v}(s; t_0).
$$

Otherwise, $a_n$ will be too large to give the correct value of $d_n$ when plugged into (70). Notice that the derivation of (69) is still valid, even though we replace (67) with (71). For instance, in Theorem 6(i) we write $v_n(t) = n^{-1/2}\sigma \times \int k(u)B(t-uh)\,du + o_P(1)$ and put $a_n = n^{-1/2}h^{1/2}, b_n = h$, so that

$$a_n^{-1}(v_n(t_0 + sb_n) - v_n(t_0))$$
$$= h^{-1/2}\sigma \int k(u)(B(t_0 + sh - uh) - B(t_0 - uh))\,du + o_P(1)$$
$$\stackrel{\mathcal{L}}{=} \sigma \int k(u)(B(s-u) - B(-u))\,du + o_P(1).$$

In the last step, we have used the fact that $B$ has stationary increments and is self similar. Thus, from (71) we obtain

$$\bar{v}(s;t_0) = \sigma \int B(u)(k(s-u) - k(-u))\,du,$$

from which follows $\beta = 1$ and

$$\tilde{v}(s;t_0) = \lim_{\delta \to 0} \sigma s \int B(u)\frac{k(s\delta - u) - k(-u)}{\delta s}\,du$$
$$= \sigma s \int B(u)k'(-u)\,du =: \sigma s Z,$$

where $Z \in N(0, \int k^2(u)\,du)$.

**6. Concluding remarks.** Several of the applications in Sections 3 and 4 are not stated in the most general form, because of a desire to keep the paper self-contained. Therefore, we would like to point out some generalizations that can be made. Section 3 on regression could be made more inclusive by allowing for heteroscedasticity and nonequidistant design points, as considered by Wright [47]; see also [33] and [34].

Furthermore, we have not made an extensive study of all possible mixing conditions, and whether, for instance, in these cases the bounds derived in Appendices B and C apply; neither have we tried to apply our results to short range dependent subordinated Gaussian data as defined in [15]. Long range dependent data limit results under exponential subordination are derived in [18]; we have not tried to apply our results to this case.

In Section 4 it is possible to prove results for estimators of a monotone density and of its derivative, by isotonization of a kernel density estimate; the calculations are similar to Section 3.3.

It is also possible to use Theorems 1 and 3 for $p \neq 2$, which would generalize the theorems in Sections 3 and 4. Since existing results, as Wright [47] and Leurgans [33], deal with independent data, this would constitute new results for weakly dependent and long range dependent data.

Unimodal density estimation, with known or unknown mode, is related to density estimation under monotonicity assumptions, and the distributional limit results are identical to ours when using the respective NPMLE. We conjecture that the results for long range dependent data in Section 4 hold also for unimodal densities.

The case $p = 1$, corresponding to nondifferentiable target functions, is not covered in this paper. It is treated in [4]; the approach is somewhat different and the limit distributions are different from the ones obtained in this paper.

It would be interesting to extend our results to process results, such as in, for example, [5, 28, 31]; we have not attempted to do so in this paper.

Alongside regression and density estimation, a third topic which can be treated with arguments similar to this paper is the monotone deconvolution problem; see [45]. In fact, in [1], the asymptotic theory of Section 2 is applied to the monotone deconvolution problem.

## APPENDIX A: PROOFS OF RESULTS IN SECTION 2

In this appendix we prove the statements in Section 2.

PROOF OF PROPOSITION 1. We first show that if Assumption A2 holds, then $g_n$ can be bounded in the following manner: for any constant $\kappa > 0$, there is a $\tau < \infty$, such that

$$(72) \qquad \liminf_{n \to \infty} \inf_{|s| \geq \tau} (g_n(s) - \kappa|s|) \geq 0.$$

To prove this, from (12) it follows that, given any $\tau > 0$ and $\varepsilon$ such that $0 < \varepsilon < A\tau^p/2$,

$$g_n(\pm\tau) \geq A\tau^p - \varepsilon,$$

if $n \geq n_0$ for some finite $n_0 = n_0(\varepsilon)$. Since $g_n(0) = 0$ and $g_n$ is convex, it follows that

$$g_n(s) \geq (A\tau^p - \varepsilon)|s|/\tau \geq \tfrac{1}{2}A\tau^{p-1}|s|,$$

when $|s| \geq \tau$, for all $n \geq n_0$. Thus, (72) holds with $\kappa = \tfrac{1}{2}A\tau^{p-1}$. Since $\tau$ can be chosen arbitrarily large, so can $\kappa$.

We are now ready to establish Assumptions A3 and A4. Choose $\tau = \tau(\delta, \varepsilon)$ as in (14). Then (14) and (72) imply

$$\liminf_{n \to \infty} P\left( \inf_{|s| \geq \tau} (y_n(s) - \kappa(1-\varepsilon)|s|) \geq 0 \right) > 1 - \delta,$$

and this proves Assumption A3 [with $\kappa(1-\varepsilon)$ in place of $\kappa$]. To establish Assumption A4, we notice that the convexity of $g_n$ and $g_n(0) = 0$ imply that $g_n(s)/s$ is increasing on $\mathbb{R}^+$. With $\tau = \tau(\delta, \varepsilon)$ as above, we get

$$\inf_{\tau \leq s \leq c} \frac{y_n(s)}{s} \leq (1+\varepsilon) \inf_{\tau \leq s \leq c} \frac{g_n(s)}{s} = (1+\varepsilon)\frac{g_n(\tau)}{\tau}$$

and, similarly,
$$\inf_{s \geq c} \frac{y_n(s)}{s} \geq (1-\varepsilon) \inf_{s \geq c} \frac{g_n(s)}{s} \geq (1-\varepsilon) \frac{g_n(c)}{c} \geq (1-\varepsilon) \frac{g_n(\tau)}{\tau}.$$

Choose $M > 0$ so that $\sup_n |g_n(\tau)| \leq M$. Then
$$\inf_{\tau \leq s \leq c} \frac{y_n(s)}{s} - \inf_{s \geq c} \frac{y_n(s)}{s} \leq 2\varepsilon \frac{g_n(\tau)}{\tau} \leq \frac{2\varepsilon M}{\tau}.$$

Since $\varepsilon$ can be made arbitrarily small, the first part of Assumption A4 is proved. The second part is proved in the same way. □

PROOF OF PROPOSITION 2. Define $h(z)$ as in (17). Then since $T(z)'(0) = T(z)'(0-)$ if and only if $h$ is a functional continuous at $z$, the set of discontinuities of $h$ is
$$D_h = \{z : T_c(z)'(0) > T_c(z)'(0-)\}$$
[recall that $T_c(z)'$ denotes the right-hand derivative of $T_c(z)$]. Let
$$\Omega_c(a, \varepsilon) = \{z : z(s) - z(0) - as \geq \varepsilon |s| \text{ for all } s \in [-c, c]\}.$$
Then $D_h \subset \bigcup_{i=1}^{\infty} \Omega_c(a_i, \varepsilon_i)$ for some countable sequence $\{(a_i, \varepsilon_i)\}$. Thus, (18) implies $P(D_h) = 0$, and Assumptions A1 and A2 imply $y_n \xrightarrow{\mathcal{L}} y$. By assumption $D_h^c$ is separable and completely regular and, thus, the continuous mapping theorem (cf. Theorem 4.12 in [38]) implies $h(y_n) \xrightarrow{\mathcal{L}} h(y)$. □

We will next go through a sequence of results that were used in the proofs of Theorems 1 and 2 and Corollary 1. We start by stating some elementary properties of the functionals. A point $t$ such that $T(y)(t) = T_c(y)(t)$ we call a point of touch of $T(y)$ and $T_c(y)$.

LEMMA A.1. *Assume $y \in D(\mathbb{R})$. If $T(y)$ and $T_c(y)$ have no points of touch on the interval $I \subset \mathbb{R}$, then $T(y)$ is linear on $I$. If $A \subset B$ are finite subsets of $\mathbb{R}$ and $y$ is bounded from above by $M$ on $A$, and bounded from below by the same $M$ on $B^c$, then*

(73) $$\inf_{s \in B} |y(s) - T(y)(s)| = 0.$$

*For any interval $O$ of $\mathbb{R}$, and functions $l, h$ on $O$ such that $l$ is linear and constant $a$ we have*

(74) $$T_O(h + l) = T_O(h) + l, \qquad T_O(ah) = aT_O(h),$$

*$T$ is monotone, that is,*

(75) $$y_1 \leq y_2 \implies T(y_1) \leq T(y_2).$$

*If $r$ is another function on $O$,*

(76) $$\sup_{t \in O} |T_O(r + h)(t) - T_O(h)(t)| \leq \sup_{t \in O} |r(t)|.$$

PROOF. Assume $T(y)$ and $T_c(y)$ have no point of touch on $I$. Then since $T_c(y)$ is a minorant of $y$, $T(y)$ and $y$ also have no point of touch on $I$ and, thus, $T(y)$ is linear on $I$.

To prove (73), suppose
$$\inf_{s \in B} |y(s) - T(y)(s)| = \varepsilon > 0.$$

Then $T(y)$ is a straight line $l$ on $B$, and further, $T(y) \leq M - \varepsilon$ on $A$. Assume w.l.o.g. that $l' \geq 0$ and let $b$ be the left end point of $B$. Then $l(s) \leq M - \varepsilon$ on $(-\infty, b]$ and, thus,
$$\inf_{s \leq b}(y(s) - l(s)) \geq \varepsilon > 0,$$

which is impossible by the construction of $T(y)$.

Equations (74) and (75) are immediate, and (76) follows from these two, since for arbitrary $t \in O$,
$$\begin{aligned}
T_O(h(t)) - \sup_{s \in O} |r(s)| &= T_O\left(h(t) - \sup_{s \in O} |r(s)|\right) \\
&\leq T_O(h(t) + r(t)) \\
&\leq T_O\left(h(t) + \sup_{s \in O} |r(s)|\right) \\
&= T_O(h(t)) + \sup_{s \in O} |r(s)|.
\end{aligned}$$
□

We next state a local limit distribution theorem.

LEMMA A.2. *Let $t_0 \in J$ be fixed, and assume Assumptions* A1 *and* A2 *hold. Then*
$$d_n^{-p}[T_{c,n}(x_n)(t_0) - x_n(t_0)] \xrightarrow{\mathcal{L}} T_c[A|s|^p + \tilde{v}(s; t_0)](0)$$
*as $n \to \infty$, and with $A > 0$ as in Assumption* A2.

PROOF. A $t$ varying in $[t_0 - cd_n, t_0 + cd_n]$ can be written as $t = t_0 + sd_n$ with $s \in [-c, c]$. Using the representation (23) and (74), we have
$$d_n^{-p}[T_{c,n}(x_n)(t_0) - x_n(t_0)] = T_c[g_n(s) + \tilde{v}_n(s; t_0)](0).$$

From Assumption A2, and the fact that $\tilde{v}_n \xrightarrow{\mathcal{L}} \tilde{v}$ on $D[-c, c]$, we get that $g_n(s) + \tilde{v}_n(s; t_0) \xrightarrow{\mathcal{L}} A|s|^p + \tilde{v}(s; t_0)$ on $D[-c, c]$. Equations (74) and the fact that $T$ is a continuous map from $D[-c, c]$ to $C[-c, c]$ [i.e., (76)] imply by the continuous mapping theorem the statement of the theorem. □

Next we show that the difference between the localized functional $T_{c,n}$ and the global $T$ goes to zero as $c$ grows to infinity.

Let us consider a sequence of stochastic processes $\{y_n\}_{n\geq 1}$ in $D(-\infty, \infty)$, for which we will state a truncation result. First we need the following additional assumption.

ASSUMPTION A.1 (Compact boundedness). For every compact set $K$ and $\delta > 0$, there is a finite $M = M(K, \delta)$ such that
$$\limsup_{n \to \infty} P\left(\sup_{s \in K} |y_n(s)| > M\right) < \delta.$$

THEOREM A.1. *Assume $y_n$ satisfies Assumptions A3, A4 and A.1. Then for every finite interval $I$ in $\mathbb{R}$ and $\varepsilon > 0$,*
$$\lim_{c \to \infty} \limsup_{n \to \infty} P\left(\sup_I |T_c(y_n)'(\cdot) - T(y_n)'(\cdot)| > \varepsilon\right) = 0,$$
$$\lim_{c \to \infty} \limsup_{n \to \infty} P\left(\sup_I |T_c(y_n)(\cdot) - T(y_n)(\cdot)| > \varepsilon\right) = 0.$$

PROOF. Let $\delta > 0$ be arbitrary and put $K = [-1, 1]$. Define the sets
$$A(n, \tau, M, \kappa) = \left\{\sup_{s \in K} |y_n(s)| < M\right\} \cap \left\{\inf_{|s| \geq \tau} (y_n(s) - \kappa|s|) > 0\right\}.$$

If $I$ is an arbitrary interval of $[-c, c]$, define the sets
$$B(n, c, I, \varepsilon) = \left\{\sup_I |T_c(y_n)'(t) - T(y_n)'(t)| < \varepsilon\right\}.$$

From Lemma A.4, it follows that

(77) $\qquad B(n, c, \{-\tau\}, \varepsilon) \cap B(n, c, \{\tau\}, \varepsilon) \subset B(n, c, I, \varepsilon)$

for any $I \subset [-\tau, \tau]$, if $\tau \leq c$.

We will show that, given any $\delta > 0$, if $c$ is large enough,

(78) $\qquad \limsup_{n \to \infty} P(B(n, c, \{\tau\}, \varepsilon)^c \cap A(n, \tau, M, \kappa)) \leq \delta,$

(79) $\qquad \limsup_{n \to \infty} P(B(n, c, \{-\tau\}, \varepsilon)^c \cap A(n, \tau, M, \kappa)) \leq \delta.$

Combining Assumptions A3 and A.1, we find that if $M$ and $\tau$ are large enough and $\kappa > 0$ is small enough, then
$$\limsup_{n \to \infty} P(A(n, \tau, M, \kappa)^c) \leq \delta.$$

Using (77), this will imply that, with a $c$ large enough and for all large enough $n$,
$$P(B(n, c, I, \varepsilon)^c) \leq P(A(n, \tau, M, \kappa)^c)$$

$$+ P(B(n,c,\{-\tau\},\varepsilon)^c \cap A(n,\tau,M,\kappa))$$
$$+ P(B(n,c,\{\tau\},\varepsilon)^c \cap A(n,\tau,M,\kappa))$$
$$\leq 2\delta + \delta + \delta,$$

and since $\delta > 0$ is arbitrary, the first assertion of the theorem will follow.

Without loss of generality, we assume that $\tau$ is chosen so large that $\tau \geq M/\kappa$. Then, given $A(n,s,M)$, we have, by our choice of $\tau$,

$$\inf_{|s|\geq \tau} y_n(s) \geq M \geq \sup_{s\in K} y_n(s).$$

Let $\zeta(\cdot;c,n,\tau)$ be the tangent line of $T_c(y_n)(s)$ at $s = \tau$, with slope $T_c(y_n)'(\tau+)$. Then exactly one of the following three events can take place. If $c > \tau$, for all large enough $n$, we have the following:

1. $\zeta(s;c,n,\tau) \leq y_n(s)$ for all $s \notin [-c,c]$;
2. $\zeta(s;c,n,\tau) \begin{cases} > y_n(s), & \text{for some } s \geq c, \\ \leq y_n(s), & \text{for all } s \leq -c; \end{cases}$
3. $\zeta(s;c,n,\tau) \begin{cases} > y_n(s), & \text{for some } s \leq -c, \\ \leq y_n(s), & \text{for all } s \geq c. \end{cases}$

In the case 1, $T_c(y_n)'(\tau) = T(y_n)'(\tau)$ if $c > \tau$.

From the assumptions defining case 2, we get, if $c > \tau$ and $A(n,\tau,M,\kappa)$ holds,

$$\begin{aligned}
\inf_{s\geq c} \frac{y_n(s) - T(y_n)(\tau)}{s - \tau} &\leq T(y_n)'(\tau) \leq T_c(y_n)'(\tau) \\
&\leq \inf_{\tau \leq s \leq c} \frac{y_n(s) - T_c(y_n)(\tau)}{s - \tau} \\
&\leq \inf_{\tau \leq s \leq c} \frac{y_n(s) - T(y_n)(\tau)}{s - \tau} \\
&\leq \left| \frac{y_n(2\tau) - T(y_n)(\tau)}{\tau} \right|,
\end{aligned} \tag{80}$$

where the last inequality holds if $c > 2\tau$. Assume that $\sup_{|s|\leq 2\tau} |y_n(s)| \leq \tilde{M}$, with $\tilde{M}$ chosen so large that this event has probability larger then $1 - \delta/4$, for all large enough $n$. Then the right-hand side of (80) is bounded by $2\tilde{M}/\tau$. Thus, $T(y_n)'(\tau) = T_c(y_n)'(\tau)$, unless

$$\inf_{s\geq c} \frac{y_n(s) - T(y_n)(\tau)}{s - \tau} \leq \frac{2\tilde{M}}{\tau}. \tag{81}$$

But if (81) holds, we get from the first half of Assumption A4 (with $\varepsilon/3$ in place of $\varepsilon$) that, with probability $\geq 1 - \delta/4$,

$$\inf_{\tilde{\tau}\leq s \leq c} \frac{y_n(s)}{s} \leq \inf_{s\geq c} \frac{y_n(s)}{s} + \frac{\varepsilon}{3} \leq 2\frac{c-\tau}{c} \inf_{s\geq c} \frac{y_n(s)}{s-\tau} + \frac{\varepsilon}{3}$$

$$\leq 2\frac{c-\tau}{c}\left(\inf_{s\geq c}\frac{y_n(s)-T(y_n)(\tau)}{s-\tau}+\frac{\tilde{M}}{c-\tau}\right)+\frac{\varepsilon}{3}$$

$$\leq 2\frac{c-\tau}{c}\left(\frac{2\tilde{M}}{\tau}+\frac{\tilde{M}}{c-\tau}\right)+\frac{\varepsilon}{3}\leq\frac{6\tilde{M}}{\tau}+\frac{\varepsilon}{3},$$

for $c$ large enough. Using Assumption A4 again with $\varepsilon/3$ in place of $\varepsilon$ and for some large $\tilde{\tau} > \tau$, this implies

$$|T_c(y_n)'(\tau) - T(y_n)'(\tau)|$$

$$\leq \inf_{\tilde{\tau}\leq s\leq c}\frac{y_n(s)-T_c(y_n)(\tau)}{s-\tau} - \inf_{s\geq\tilde{\tau}}\frac{y_n(s)-T(y_n)(\tau)}{s-\tau}$$

$$\leq \frac{2\tilde{M}}{\tilde{\tau}-\tau} + \inf_{\tilde{\tau}\leq s\leq c}\frac{y_n(s)}{s-\tau} - \inf_{s\geq\tilde{\tau}}\frac{y_n(s)}{s-\tau}$$

$$\leq \frac{2\tilde{M}}{\tilde{\tau}-\tau} + \left(\inf_{\tilde{\tau}\leq s\leq c}\frac{y_n(s)}{s} - \inf_{s\geq\tilde{\tau}}\frac{y_n(s)}{s}\right) + \left(\inf_{\tilde{\tau}\leq s\leq c}\frac{y_n(s)}{s-\tau} - \inf_{\tilde{\tau}\leq s\leq c}\frac{y_n(s)}{s}\right)$$

$$\leq \frac{2\tilde{M}}{\tilde{\tau}-\tau} + \left(\inf_{\tilde{\tau}\leq s\leq c}\frac{y_n(s)}{s} - \inf_{s\geq\tilde{\tau}}\frac{y_n(s)}{s}\right) + \left(\frac{\tilde{\tau}}{\tilde{\tau}-\tau}-1\right)\inf_{\tilde{\tau}\leq s\leq c}\frac{y_n(s)}{s}$$

$$\leq \frac{2\tilde{M}}{\tilde{\tau}-\tau} + \frac{\varepsilon}{3} + \frac{\tau}{\tilde{\tau}-\tau}\left(\frac{6\tilde{M}}{\tau}+\frac{\varepsilon}{3}\right) \leq \varepsilon,$$

provided $\tilde{\tau}$ is chosen large enough and then $c$ is chosen so large that Assumption A4 holds. Thus, given $A(n,\tau,M,\kappa)$ and case 2, $B(n,c,\{\tau\},\varepsilon)$ holds unless $\sup_{|s|\leq 2\tau}|y_n(s)| > \tilde{M}$ or if the first half of Assumption A4 fails, which is an event with probability at most $\delta/4 + \delta/4 = \delta/2$.

Given $A(n,\tau,M,\kappa)$ and case 3, a similar argument implies that $B(n,c,\{\tau\},\varepsilon)$ fails with probability at most $\delta/2$. Combining cases 1–3, we deduce (78), and (79) is proved in the same way.

To show the second part of the theorem, for $I$ an arbitrary interval of $(-c,c)$ containing $[-\tau,\tau]$, define the sets

$$C(n,c,I,\varepsilon) = \left\{\sup_{s\in I}|T_c(y_n)(s) - T(y_n)(s)| < \varepsilon\right\}.$$

Let $L = \text{length}(I)$. Suppose that Assumptions A3 and A.1 hold (an event with probability $\geq 1 - 2\delta$). We will apply (73) with $A = K$ as in Assumption A.1 and $B = [-\tau,\tau]$. Assume also $\tau \geq M/\kappa$, with $M,\tau,\kappa$ as in Assumptions A3 and A.1. Then there is an $\eta \in [-\tau,\tau]$ such that $|y_n(\eta) - T(y_n)(\eta)| \leq \varepsilon/2$. Since $T(y_n)(\eta) \leq T_c(y_n)(\eta) \leq y_n(\eta)$, we get $|T_c(y_n)(\eta) - T(y_n)(\eta)| \leq \varepsilon/2$. Thus,

$$P(C(n,c,I,\varepsilon)) \leq P(B(n,c,I,\varepsilon/2L)) + 2\delta,$$

and the second part of the theorem follows from the first. □

NOTE A.1. If $T$ is replaced with the greatest convex minorant $T_O$ on an interval $O$ on $\mathbb{R}$, Theorem A.1 trivially still holds. If $T$ is replaced with $T_{O_n}$ where $O_n$ is a sequence of intervals such that $O_n \uparrow \mathbb{R}$, the theorem still holds. In the latter case Assumptions A3, A4 and A.1 may be relaxed somewhat; all suprema and infima with respect to $s$ over $\mathbb{R}$ can be relaxed to suprema and infima with respect to $s$ over $O_n$.

LEMMA A.3. *Assume Assumptions* A1, A2, A3 *and* A4 *hold. Define* $A_{n,\Delta} = [t_0 - \Delta d_n, t_0 + \Delta d_n]$. *Then for every finite* $\Delta$, *and every* $\varepsilon > 0$,

$$\lim_{c \to \infty} \liminf_{n \to \infty} \mathbf{P}\left[\sup_{A_{n,\Delta}} d_n^{-p}|T_{c,n}(x_n)(\cdot) - T_J(x_n)(\cdot)| \leq \varepsilon\right] = 1,$$

$$\lim_{c \to \infty} \liminf_{n \to \infty} \mathbf{P}\left[\sup_{A_{n,\Delta}} d_n^{-p+1}|T_{c,n}(x_n)'(\cdot) - T_J(x_n)'(\cdot)| \leq \varepsilon\right] = 1.$$

PROOF. From (74) and (23), it follows that

$$\sup_{A_{n,\Delta}} d_n^{-p}|T_{c,n}(x_n)(\cdot) - T_J(x_n)(\cdot)| = \sup_{[-\Delta,\Delta]} |T_c(y_n)(\cdot) - T_{J_{n,t_0}}(y_n)(\cdot)|$$

and

$$\sup_{A_{n,\Delta}} d_n^{-p+1}|T_{c,n}(x_n)'(\cdot) - T_J(x_n)'(\cdot)| = \sup_{[-\Delta,\Delta]} |T_c(y_n)'(\cdot) - T_{J_{n,t_0}}(y_n)'(\cdot)|,$$

with $y_n$ as defined in (13).

If $J = \mathbb{R}$, we use Theorem A.1 with $I = [-\Delta, \Delta]$, and if $J \neq \mathbb{R}$, we use Note A.1 with $O_n = J_{n,t_0}$.

Assumptions A3 and A4 are satisfied because of Proposition 1, and Assumption A.1 is implied by Assumption A1 and (12). Thus, all the regularity conditions of Theorem A.1 are satisfied, and the lemma follows. □

LEMMA A.4. *Suppose* $y \in D(\mathbb{R})$. *Then the function* $t \mapsto T_c(y)'(t) - T(y)'(t)$ *is increasing on* $[-c, c]$.

PROOF. Let $\{J_k\}$ be a sequence of open intervals in $(-c, c)$ such that their union covers $(-c, c)$. Without loss of generality, we can assume that each $J_k$ either contains no points of touch of $T(y)$ and $T_c(y)$ or it contains exactly one simply connected set $\Omega_{J_k}$ of points of touch (so then $\Omega_{J_k}$ consists of either a simple point or it is an interval). If $\Omega_{J_k}$ is empty, by the first part of Lemma A.1, $T(y)$ is linear on $J_k$, and since $T_c(y)$ is convex, the assertion follows. If $\Omega_{J_k}$ is nonempty, since $T_c(y) \geq T(y)$, for all $t \in [-c, c]$, the assertion holds again, and the lemma is proven. □

TABLE 2
*Rescaled processes*

| Section | $\tilde{v}_n(s;t)$ | $c_n$ | $\hat{n}$ |
|---|---|---|---|
| 3.1 | $c_n w_{\hat{n}}(s)$ | $d_n^{-2}(nh)^{-1}\sigma_{\hat{n}}$ | $nd_n$ |
| 3.2.1 | $c_n \int (\bar{w}_{\hat{n}}(s-u) - \bar{w}_{\hat{n}}(-u))k'(-u)\,du$ | $d_n^{-2}(nh)^{-1}\sigma_{\hat{n}}$ | $nh$ |
| 3.3.1 | $c_n \int (w_{\hat{n}}(s-u) - w_{\hat{n}}(-u))k(u)\,du$ | $d_n^{-2}(nh)^{-1}\sigma_{\hat{n}}$ | $nh$ |
| 3.3.2 | $c_n \int (w_{\hat{n}}(-uh/d_n + s) - w_{\hat{n}}(-uh/d_n)k(u))\,du$ | $d_n^{-1}n^{-1}\sigma_{\hat{n}}$ | $nd_n$ |
| 3.3.3 | $c_n \int w_{\hat{n}}(u)\frac{k(sd_n/h-u)-k(-u)}{d_n/h}\,du$ | $(d_n nh)^{-1}\sigma_{\hat{n}}$ | $nh$ |

## APPENDIX B: BOUND ON DRIFT OF PROCESS PART: PARTIAL SUMS

In this appendix we will establish Proposition 1 for the various applications of Section 3. By Proposition 1, Assumptions A3 and A4 are implied by Assumption A2 and the following:

ASSUMPTION B.1. Assume that for any $\varepsilon, \delta > 0$, there exist $\kappa = \kappa(\varepsilon, \delta) > 0$ and $\tau = \tau(\varepsilon, \delta) > 0$ such that

$$\sup_n P\bigg(\sup_{|s|\geq \tau} \frac{\tilde{v}_n(s;t)}{\kappa|s|} > \varepsilon\bigg) < \delta.$$

In all the cases of Section 3, the rescaled process $\tilde{v}_n$ can be written as a function of the partial sum process $w_{\hat{n}}$; see Table 2. (Note that in some of the cases in Section 3, $\tilde{v}_n$ is a function of $\bar{w}_{\hat{n}}$ instead of $w_{\hat{n}}$. However, since $\bar{w}_{\hat{n}}$ is a smoothed version of $w_{\hat{n}}$, the bounds on $w_{\hat{n}}$ established in this section are easily shown to translate to bounds on $\bar{w}_{\hat{n}}$.)

In all the above cases we have $c_n \to c > 0$ as $n \to \infty$. Therefore, we start by establishing the following:

LEMMA B.1. *Suppose $\{\varepsilon_i\}_{i\geq 1}$ is a stationary independent process, a weakly dependent sequence satisfying Assumption A8 or A9 or a long range dependent subordinated Gaussian sequence with parameters d and r, and $\beta$ as in (36), and assume that $E(\varepsilon_i) = 0$ and $\text{Var}(\varepsilon_i) = \sigma^2 < \infty$. Then for each $\varepsilon, \delta, \kappa > 0$, there exist $\tau = \tau(\varepsilon, \delta, \kappa) > 0$ and $m_0 = m_0(\varepsilon, \delta, \kappa) < \infty$ such that*

(82) $$\sup_{m\geq m_0} P\bigg(\sup_{|s|\geq \tau} \frac{|w_{\hat{n}}(s)|}{\kappa|s|} > \varepsilon\bigg) < \delta.$$

Observe that, with $\tau < \infty, \varepsilon > 0$ and $\hat{n}$ fixed, we have

$$P\bigg(\sup_{|s|\geq \tau} \frac{|w_{\hat{n}}(s)|}{\kappa|s|} > \varepsilon\bigg)$$

$$
\begin{aligned}
&\leq 2P\bigg(\sup_{s\geq\tau}\frac{|w_{\hat{n}}(s)|}{\kappa|s|} > \varepsilon\bigg) \\
&\leq 2\sum_{s_i\geq\tau} P\bigg(|w_{\hat{n}}(s_i)| > \frac{\varepsilon}{2}\kappa s_i\bigg) \\
&\quad + 2\sum_{s_i\geq\tau} P\bigg(\sup_{s_{i-1}\leq s\leq s_i}|w_{\hat{n}}(s) - w_{\hat{n}}(s_{i-1})| > \frac{\varepsilon}{2}\kappa s_{i-1}\bigg),
\end{aligned}
\tag{83}
$$

where $\cdots < s_{-1} < 0 < s_1 < \cdots$ is a partition of $\mathbb{R}$. Note that we used the fact that $s \mapsto \kappa|s|$ is increasing for $s > 0$ in the second inequality in (83).

LEMMA B.2. *Under the assumptions in Lemma* B.1,

$$P\bigg(\sup_{s_i\leq s\leq s_{i+1}}|w_{\hat{n}}(s) - w_{\hat{n}}(s_i)| > \frac{\varepsilon}{2}\kappa|s_i|\bigg) \leq \frac{C\Delta_i^{2\beta}}{\varepsilon^2 s_i^2},$$

where $\Delta_i = s_{i+1} - s_i$, and $\beta = 1/2$ in the independent and weakly dependent cases and $1/2 < \beta = 1 - rd/2 < 1$ in the long range dependent case.

PROOF. Let $S_k = \sum_{i=1}^k \varepsilon_i$, $\tilde{n} = \hat{n}\Delta_i$ and assume that $s_i \leq s \leq s_{i+1}$. Then by stationarity,

$$w_{\hat{n}}(s) - w_{\hat{n}}(s_i) \stackrel{\mathcal{L}}{=} w_{\hat{n}}(s - s_i),$$

at least when $s_i = d_n^{-1}(t_i - t_0)$ for some observation point $t_i$, which we assume w.l.o.g. Thus,

$$P\bigg(\sup_{s_i\leq s\leq s_{i+1}}|w_{\hat{n}}(s) - w_{\hat{n}}(s_{i+1})| > \frac{\varepsilon}{2}\kappa|s_i|\bigg) = P\bigg(\max_{1\leq k\leq \tilde{n}}|S_k| > \lambda\sigma_{\tilde{n}}\bigg),$$

with $\sigma_{\tilde{n}}^2 = \mathrm{Var}(S_{\tilde{n}})$ and

$$\lambda = \frac{\varepsilon}{2}\frac{\sigma_{\hat{n}}}{\sigma_{\tilde{n}}}\kappa|s_i|.$$

For independent data, equation (10.7) in [8] implies

$$P\bigg(\max_{1\leq k\leq \tilde{n}}|S_k| > \lambda\sigma_{\tilde{n}}\bigg) \leq \frac{C}{\lambda^2}. \tag{84}$$

Since $\sigma_{\tilde{n}}/\sigma_m = \Delta_i^{1/2}$ for independent data, this proves the lemma in this case.

In the weakly dependent mixing case, we use the results of McLeish [36] to prove (84). Thus, denoting $\|z\|_q = (E|z|^q)^{1/q}$ for a random variable $z$, we call the sequence $\{\varepsilon_i\}$ a mixingale if

$$\|\varepsilon_n - E(\varepsilon_n|\mathcal{F}_{n+\hat{n}})\|_2 \leq \psi_{\hat{n}+1}c_n, \tag{85}$$

$$\|E(\varepsilon_n|\mathcal{F}_{n-\hat{n}})\|_2 \leq \psi_{\hat{n}}c_n, \tag{86}$$

ORDER RESTRICTED INFERENCE 51

with $c_n, \psi_{\hat{n}}$ finite and nonnegative constants and $\lim_{\hat{n}\to\infty} \psi_{\hat{n}} = 0$. Since each $\varepsilon_n$ is $\mathcal{F}_n$-measurable, (85) holds. Assuming that $\varepsilon_i$ has finite variance, or in the case of $\alpha$-mixing, finite fourth moment, Lemma 2.1 in [36] implies

$$\|E(\varepsilon_n|\mathcal{F}_{n-\hat{n}})\|_2 \le 2\phi(\hat{n})^{1/2}\|\varepsilon_n\|_2,$$
$$\|E(\varepsilon_n|\mathcal{F}_{n-\hat{n}})\|_2 \le 2(\sqrt{2}+1)\alpha(\hat{n})^{1/4}\|\varepsilon_n\|_4.$$

Using Assumption A8 or A9, we will apply Theorem 1.6 in [36] with either $\psi_{\hat{n}} = 2\phi(\hat{n})^{1/2}$ and $c_n = \|\varepsilon_n\|_2 = \sigma$ or $\psi_{\hat{n}} = 2(\sqrt{2}+1)\alpha(\hat{n})^{1/4}$ and $c_n = \|\varepsilon_n\|_4$. In either case $\sum_{\hat{n}=1}^{\infty} \psi_{\hat{n}}^{2-\varepsilon} < \infty$ for some $\varepsilon > 0$, which shows that $\psi(\hat{n})$ is of size $-1/2$ with the McLeish [36] terminology (as noted on top of page 831 in [36]). Notice also that $\sum_{i=1}^{\hat{n}} c_i^2 = C\hat{n}$ for some $C > 0$ and $\sigma_{\hat{n}}^2 \sim \kappa^2 \hat{n}$ according to (33). Thus, Theorem 6.1 in [36] and Chebyshev's inequality imply (84). Notice that $\sigma_{\tilde{n}}/\sigma_{\hat{n}} \sim \Delta_i^{1/2}$ as $\hat{n}, \tilde{n} \to \infty$.

In the long range dependent case we use Theorem 12.2 in [8]. Thus,

$$E(S_{\tilde{n}}^2) \sim \eta_r^2 l_1(\tilde{n}) \tilde{n}^{2\beta},$$

with $l_1$ as in (35), and according to de Haan [16], equation (12.42) in [8] is satisfied, with

$$\gamma = 2,$$
$$\alpha = 2\beta,$$
$$u_l = (C_1 \eta_r^2 l_1(\tilde{n}))^{1/2\beta},$$

for some constant $C_1 > 0$. Theorem 12.2 in [8] then implies that

$$(87) \qquad P\left(\max_{1\le k\le \tilde{n}} |S_k| > \lambda\sigma_{\tilde{n}}\right) \le \frac{K'_{2,2\beta}}{(\lambda\sigma_{\tilde{n}})^2}\left(\sum_{i=1}^{\tilde{n}} u_i\right)^{2\beta} = \frac{C}{\lambda^2},$$

with $C = C_1 K'_{2,2\beta}$. From de Haan [16] it follows that

$$C_1^{-1}\Delta_i^{2\beta} \le \frac{\sigma_{\tilde{n}}^2}{\sigma_{\hat{n}}^2} = \frac{l_1(\hat{n}\Delta_i)(\hat{n}\Delta_i)^{2\beta}}{l_1(\hat{n})\hat{n}^{2\beta}} \le C_1\Delta_i^{2\beta},$$

for all large enough $\hat{n}$. Substituting for $\lambda$ in (87) completes the proof. □

PROOF OF LEMMA B.1. From (83) and Lemma B.2, we have

$$P\left(\sup_{|s|\ge\tau} \frac{|w_{\hat{n}}(s)|}{\kappa|s|} > \varepsilon\right) \le 8\sum_{s_i>\tau} \frac{Ew_{\hat{n}}^2(s_i)}{\varepsilon^2\kappa^2 s_i^2} + 2\sum_{s_{i+1}>\tau} \frac{C\Delta_i^{2\beta}}{\varepsilon^2 s_i^2}.$$

In order to take care of the slowly varying factor in $Ew_{\hat{n}}^2(s_i)$ for long range dependent data, we write

$$Ew_{\hat{n}}^2(s_i) \le Cs_i^{2\beta'},$$

where $\beta' = 1/2$ in the independent and weakly dependent case and $\beta < \beta' < 1$ in the long range dependent case (this is no restriction if we assume $\Delta_i \geq \delta > 0$ for all $i$ and some constant $\delta$). Replacing $\Delta_i^{2\beta}$ with $\Delta_i^{2\beta'}$, we want to examine whether the sums

$$\sum_{s_i > \tau} \frac{s_i^{2\beta'}}{s_i^2} = \sum_{s_i > \tau} s_i^{2\beta'-2},$$

$$\sum_{s_{i+1} > \tau} \frac{\Delta_i^{2\beta'}}{s_i^2}$$

tend to zero as $\tau \to \infty$. Clearly, choosing $\Delta_i = \delta$, the first sum is divergent. Instead we let $s_i = i^\rho$ with $\rho > 1$. Thus, the sums are of the order

$$\sum_{i > \tau^{1/\rho}} i^{\rho(2\beta'-2)} \sim (\tau^{1/\rho})^{2\rho(\beta'-1)+1},$$

$$\sum_{i > \tau^{1/\rho}-1} i^{2(\rho-1)\beta'-2\rho} \sim (\tau^{1/\rho})^{2(\rho-1)\beta'-2\rho+1}.$$

The demand that both expressions should converge to zero as $\tau \to \infty$ implies that

$$-2(1-\beta') + \frac{1}{\rho} < 0,$$

$$-2(1-\beta') - \frac{2\beta'-1}{\rho} < 0,$$

which shows that we should choose $\rho > 1/2(1-\beta')$, and this completes the proof. □

Lemma B.1 immediately proves that $\tilde{v}_n(\cdot; t)$ in Section 3.1 satisfies Assumption B.1.

LEMMA B.3. *Assume $\{w_n\}$ satisfies (82) in Lemma B.1 and is uniformly bounded in probability on compact intervals. Let $\{l_n\}$ be a sequence of functions with $\mathrm{supp}(l_n) \subset [-K, K]$ for some $K > 0$ and all $n$, and with $\sup_n \int |l_n(u)|\, du < \infty$. Then*

$$\tilde{v}_n(s) = \int (w_{\hat{n}}(s-u) - w_{\hat{n}}(-u)) l_n(u)\, du$$

*satisfies Assumption B.1.*

PROOF. Since

$$|\tilde{v}_n(s)| \leq \left( \sup_{u[-K,K]} |w_m(s-u)| + \sup_{u \in [-K,K]} |w_m(-u)| \right) \int |l_n(u)|\, du,$$

we obtain

$$\sup_{|s|\geq\tau} \frac{|\tilde{v}_n(s)|}{\kappa|s|} \leq C\bigg(\sup_{|s|\geq\tau}\sup_{u\in[-K,K]}\frac{|w_m(s-u)|}{\kappa|s|}$$
$$+ \sup_{|s|\geq\tau}\sup_{u\in[-K,K]}\frac{|w_m(-u)|}{\kappa|s|}\bigg)$$
$$\leq C'\sup_{|s|\geq\tau}\sup_{u\in[-K,K]}\frac{|w_m(s-u)|}{\kappa|s-u|} + C\sup_{u\in[-K,K]}\frac{|w_m(-u)|}{\kappa|\tau|}$$
$$= C'\sup_{|s|\geq\tau-K}\frac{|w_m(s)|}{\kappa|s|} + C\sup_{u\in[-K,K]}\frac{|w_m(-u)|}{\kappa|\tau|},$$

with $C = \sup_n \int |l_n(u)|\,du$ and $C' = C\sup_{|s|\geq\tau}\sup_{u\in[-K,K]}\frac{|s-u|}{|s|}$. Finally, (82) and the fact that $\{w_m\}$ is uniformly bounded in probability on compact intervals finish the proof. □

Applying Lemma B.3 with $l_n(u)$ equal to $c_n k'(-u), c_n k(u)$ and $c_n k(d_n u/h) \times d_n/h$, respectively, establishes Assumption B.1 in Sections 3.2.1, 3.3.1 and 3.3.2.

LEMMA B.4. *Assume $\{w_n\}$ satisfies* (82) *in Lemma* B.1. *Let $l$ be a function of bounded variation with support in $[-1,1]$, and assume $\{\rho_n\}$ is a sequence of numbers such that $\lim_{n\to\infty}\rho_n = 0$. Then*

$$\tilde{v}_n(s) = \frac{w_m * l(s\rho_n) - w_m * l(0)}{\rho_n}$$
$$= \int w_m(-u)\frac{l(s\rho_n + u) - l(u)}{\rho_n}\,du$$

*satisfies Assumption* B.1.

PROOF. We will give different bounds on $\tilde{v}_n(s)$ for small and large values of $|s|$. Assume that $|s| \leq (\tau+1)\rho_n^{-1}$, where $\tau > 0$ is a constant that will be chosen below. Then

$$|v_n(s)| \leq |s|\sup_{|u|\leq\tau+2}|w_m(u)|\int\frac{|l(s\rho_n - u) - l(-u)|}{|s|\rho_n}\,du$$
$$\leq |s|\sup_{|u|\leq\tau+2}|w_m(u)|\int|l'(u)|\,du.$$

If instead $|s| > (\tau+1)\rho_n^{-1}$, choose arbitrary $\varepsilon, \delta, \kappa_0 > 0$. Then Lemma B.2 implies the existence of $\tau = \tau(\varepsilon, \delta, \kappa_0) > 0$ such that, with probability larger than $1 - \delta$, we have

$$|v_n(s)| \leq \rho_n^{-1}\left(\int |w_m(u)l(-u)|\,du + \varepsilon\kappa_0\int |u||l(s\rho_n - u)|\,du\right)$$

$$\leq \rho_n^{-1}\int |l(u)|\,du\left(\sup_{|u|\leq 1}|w_m(u)| + \varepsilon\kappa_0(|s|\rho_n + 1)\right)$$

$$\leq |s|\int |l(u)|\,du\left(\frac{\sup_{|u|\leq 1}|w_m(u)|}{\tau + 1} + \varepsilon\kappa_0\left(1 + \frac{1}{1+\tau}\right)\right).$$

Thus, with probability larger than $1 - \delta$,

$$\sup_{s\neq 0}\frac{|v_n(s)|}{\kappa_0|s|} \leq \max\left(\sup_{|u|\leq \tau+2}|w_m(u)|\int |l'(u)|\,du,\right.$$

$$\left.\int |l(u)|\,du\left(\frac{\sup_{|u|\leq 1}|w_m(u)|}{\tau + 1} + \varepsilon\kappa_0\frac{\tau + 2}{\tau + 1}\right)\right).$$

Since we assume that $w_m$ is bounded on compacta uniformly over $m$, with probability larger than $1 - 2\delta$, the right-hand side is bounded from above by a constant $C = C(\varepsilon, \delta, \kappa_0) > 0$. Pick $\kappa = \kappa_0 C/\varepsilon$. Then

$$\sup_{|s|\neq 0}\frac{|v_n(s)|}{\kappa|s|} \leq \frac{\varepsilon}{C}\sup_{s\neq 0}\frac{v_n(s)}{\kappa_0|s|} \leq \varepsilon,$$

with a probability larger than $1 - 2\delta$. □

Applying Lemma B.4 with $\rho_n = d_n/h$ and $l(n)$ equal to $k(u)$ establishes Assumption B.1 in Section 3.3.3.

## APPENDIX C: BOUND ON DRIFT OF PROCESS PART: EMPIRICAL DISTRIBUTIONS

In this appendix we will establish Assumption B.1 for the various applications treated in Section 4. The processes $\tilde{v}_n(s; t_0)$ are functions of $w_{n,\delta_n}(s) := w_{n,\delta_n}(s; t_0)$ for all cases treated in Section 4, as seen in Table 3.

In all the above cases we have $c_n \to c > 0$ as $n \to \infty$.

TABLE 3
*Rescaled processes*

| Section | $\tilde{v}_n(s; t)$ | $c_n$ | $\delta_n$ |
|---|---|---|---|
| 4.1 | $c_n w_{n,\delta_n}(s)$ | $d_n^{-2}(nh)^{-1}\sigma_{n,\delta_n}$ | $d_n$ |
| 4.2.1 | $c_n \int (w_{n,\delta_n}(s - u) - w_{n,\delta_n}(-u))k'(-u)\,du$ | $d_n^{-2}(nh)^{-1}\sigma_{n,\delta_n}$ | $h$ |

LEMMA C.1. *Assume $\{t_i\}$ is a stationary sequence with marginal distribution $F$, such that $f(t_0)$ exists, and $\delta_n \downarrow 0, n\delta_n \uparrow \infty$ as $n \to \infty$. Then, if $\{t_i\}$ is an independent or $\phi$-mixing sequence with $\sum_{i=1}^{\infty} n\phi(n)^{1/2} < \infty$, there exists for each $\varepsilon, \delta, \kappa > 0$ a $\tau = \tau(\varepsilon, \delta, \kappa) > 0$ such that*

(88) $$P\left(\sup_{|s|\geq\tau} \frac{|w_{n,\delta_n}(s;t_0)|}{\kappa|s|} > \varepsilon\right) < \delta,$$

*for all large enough $n$. If $\{t_i\}$ is a long range dependent subordinated Gaussian sequence and satisfies the assumptions of Theorem 9, then for each $\varepsilon, \delta > 0$, there exist $\kappa = \kappa(\varepsilon, \delta) > 0$ and $\tau = \tau(\varepsilon, \delta) > 0$ such that (88) holds.*

PROOF. We start by proving the lemma for long range dependent data. Then

$$w_{n,\delta_n}(s;t_0) = \frac{\eta_r(t_0 + s\delta_n) - \eta_r(t_0)}{\delta_n |\eta'_r(t_0)|}(1 + o(1))\sigma_n^{-1} \sum_{i=1}^{n} h_r(\xi_i)$$
$$+ C_n \delta_n^{-1}(S_n(t_0 + s\delta_n) - S_n(t_0)),$$

where $C_n \to r!/|\eta_r(t_0)|$ as $n \to \infty$. Clearly, $\|\eta_r\|_\infty = \sup_t |\eta_r(t)| < \infty$. Moreover, since $\eta'_r(t_0) \neq 0$, there exists a $\tilde{\delta} > 0$ such that $|\eta_r(t_0 + s) - \eta_r(t_0)|/|s| \leq 2|\eta'_r(t_0)|$ whenever $|s| \leq \tilde{\delta}$. Thus,

$$\frac{|\eta_r(t_0 + s\delta_n) - \eta_r(t_0)|}{\delta_n |\eta'_r(t_0)|} \leq \max\left(2, \frac{2\|\eta'_r\|_\infty}{\tilde{\delta}|\eta'_r(t_0)|}\right) \cdot |s|.$$

Further, since

$$\sigma_n^{-1} \sum_{i=1}^{n} h_r(\xi_i) \xrightarrow{\mathcal{L}} z_{r,\beta}$$

as $n \to \infty$, and since from [12],

$$\sup_s \delta_n^{-1}|S_n(t_0 + s\delta_n) - S_n(t_0)| \xrightarrow{\text{a.s.}} 0$$

as $n \to \infty$, the result follows.

In the independent and weakly dependent data case, we consider w.l.o.g. the supremum for $s \geq \tau$ only. Analogously to (83), we have

$$P\left(\sup_{s\geq\tau} \frac{|w_{n,\delta_n}(s)|}{\kappa|s|} > \varepsilon\right)$$
$$\leq \sum_{s_i \geq \tau} P\left(|w_{n,\delta_n}(s_i)| > \frac{\varepsilon}{2}\kappa s_i\right)$$
$$+ \sum_{s_i \geq \tau} P\left(\sup_{s_{i-1}\leq s\leq s_i} |w_{n,\delta_n}(s) - w_{n,\delta_n}(s_{i-1})| > \frac{\varepsilon}{2}\kappa s_{i-1}\right),$$

where $0 < s_1 < s_2 < \cdots$ is an increasing sequence. Assume first that $F \sim U(0,1)$. To proceed further, we need the following lemma, proved in [3].

LEMMA C.2. *Suppose $\{t_i\}$ is an independent or weakly dependent sequence of random variables, satisfying the assumptions of Lemma C.1. Then*

$$P\bigg(\sup_{s_i \leq s \leq s_{i+1}} |w_{n,\delta_n}(s) - w_{n,\delta_n}(s_i)| \geq \lambda\bigg) \leq K\Delta_i^2\bigg(\frac{1}{\lambda^4} + \frac{1}{\lambda^5}\bigg)$$

*for all $\lambda > 0$ if $\Delta_i = s_{i+1} - s_i \geq 1$, with $K$ a constant depending only on $\{\phi_n\}$.*

The next part of the proof proceeds similarly to the proof of Lemma B.2, so we highlight only the differences. Let $g_s$ be defined as in the proof of Lemma B.1. Then from [8], page 172, we get

$$E(w_{n,\delta_n}^2) \leq \sigma_{n,\delta_n}^{-2}\bigg(1 + 4\sum_{i=1}^{\infty} \phi_i^{1/2}\bigg) nE(g_s(t_1) - g_0(t_1))^2$$

$$\leq 2\bigg(1 + 4\sum_{i=1}^{\infty} \phi_i^{1/2}\bigg) s =: Cs.$$

By Lemma C.2 and Chebyshev's inequality, the lemma is proved for $F \sim U(0,1)$ if we can prove that the sums

$$\sum_{s_i \geq \tau} \frac{s_i}{s_i^2}, \sum_{s_i^2 \geq \tau} \frac{\Delta_i^2}{s_i^4} \text{ and } \sum_{s_i^2 \geq \tau} \frac{\Delta_i^2}{s_i^5}$$

tend to zero as $\tau \to \infty$. But this is true if $s_i = i^\rho$ for any $\rho > 1$.

Consider again a general $F$ with $f(t_0) > 0$. Let $w_{n,\delta_n}^U$ and $\sigma_{n,\delta_n}^U$ be the quantities corresponding to $w_{n,\delta_n}$ and $\sigma_{n,\delta_n}$ when $F \sim U(0,1)$. Then

(89) $$w_{n,\delta_n}(s;t_0) = w_{n,\hat{\delta}_n}^U(\hat{s}; F(t_0)),$$

where $\hat{\delta}_n = F(t_0 + \delta_n) - F(t_0)$ and $\hat{s} = (F(t_0 + s\delta_n) - F(t_0))/(F(t_0 + \delta_n) - F(t_0))$. Choose $\hat{\delta} > 0$ such that $f(t_0)/2 < |F(t_0 + s\delta_n) - F(t_0)|/|s| < 2f(t_0)$ if $|s| \leq \hat{\delta}$. Then, since $0 \leq F(t_0 + s\delta_n) - F(t_0) \leq 1$, it follows that

(90) $$\sup_{s \neq 0} \frac{\hat{s}}{s} \leq \max\bigg(4, \frac{2}{\hat{\delta} f(t_0)}\bigg)$$

for all $n$ so large that $\delta_n \leq \hat{\delta}$. Now (89), (90) and the proof of (88) when $F \sim U(0,1)$ finish the proof of (88) for general $F$. □

To establish Assumption B.1 for the various choices of $\tilde{v}_n(\cdot; t_0)$ in the table in this appendix, we proceed as in Appendix B, making use of Lemmas B.3 and C.1.

**Acknowledgments.** We would like to thank an Associate Editor and referee for their valuable comments that significantly improved the readability of the paper.

## REFERENCES

[1] ANEVSKI, D. (1999). Deconvolution under monotonicity assumptions. Technical Report 17, Centre for Mathematical Sciences, Lund Univ.

[2] ANEVSKI, D. (2003). Estimating the derivative of a convex density. *Statist. Neerlandica* **57** 245–257. MR2028914

[3] ANEVSKI, D. and HÖSSJER, O. (2000). A general asymptotic scheme for inference under order restrictions. Technical Report 1, Centre for Mathematical Sciences, Lund Univ.

[4] ANEVSKI, D. and HÖSSJER, O. (2002). Monotone regression and density function estimation at a point of discontinuity. *J. Nonparametr. Statist.* **14** 279–294. MR1905752

[5] BANERJEE, M. and WELLNER, J. A. (2001). Likelihood ratio tests for monotone functions. *Ann. Statist.* **29** 1699–1731. MR1891743

[6] BERAN, J. (1992). Statistical methods for data with long-range dependence (with discussion). *Statist. Sci.* **7** 404–427.

[7] BICKEL, P. J. and FAN, J. (1996). Some problems on the estimation of unimodal densities. *Statist. Sinica* **6** 23–45. MR1379047

[8] BILLINGSLEY, P. (1968). *Convergence of Probability Measures*. Wiley, New York. MR0233396

[9] BRADLEY, R. C. (1986). Basic properties of strong mixing conditions. In *Dependence in Probability and Statistics* (E. Eberlein and M. S. Taqqu, eds.) 165–192. Birkhäuser, Boston. MR0899990

[10] BRUNK, H. D. (1958). On the estimation of parameters restricted by inequalities. *Ann. Math. Statist.* **29** 437–454. MR0132632

[11] BRUNK, H. D. (1970). Estimation of isotonic regression. In *Nonparametric Techniques in Statistical Inference* (M. L. Puri, ed.) 177–197. Cambridge Univ. Press, London. MR0277070

[12] CSÖRGŐ, S. and MIELNICZUK, J. (1995). Density estimation under long-range dependence. *Ann. Statist.* **23** 990–999. MR1345210

[13] CSÖRGŐ, S. and MIELNICZUK, J. (1995). Distant long-range dependent sums and regression estimation. *Stochastic Process. Appl.* **59** 143–155. MR1350260

[14] CSÖRGŐ, S. and MIELNICZUK, J. (1995). Nonparametric regression under long-range dependent normal errors. *Ann. Statist.* **23** 1000–1014. MR1345211

[15] CSÖRGŐ, S. and MIELNICZUK, J. (1996). The empirical process of a short-range dependent stationary sequence under Gaussian subordination. *Probab. Theory Related Fields* **104** 15–25. MR1367664

[16] DE HAAN, L. (1970). *On Regular Variation and Its Application to the Weak Convergence of Sample Extremes*. Math. Centrum, Amsterdam. MR0286156

[17] DEHLING, H. and TAQQU, M. S. (1989). The empirical process of some long-range dependent sequences with an application to $U$-statistics. *Ann. Statist.* **17** 1767–1783. MR1026312

[18] GAJEK, L. and MIELNICZUK, J. (1999). Long- and short-range dependent sequences under exponential subordination. *Statist. Probab. Lett.* **43** 113–121. MR1693261

[19] GASSER, T. and MÜLLER, H.-G. (1984). Estimating regression functions and their derivatives by the kernel method. *Scand. J. Statist.* **11** 171–185. MR0767241

[20] GRENANDER, U. (1956). On the theory of mortality measurement. II. *Skand. Aktuarietidskr.* **39** 125–153. MR0093415
[21] GROENEBOOM, P. (1985). Estimating a monotone density. In *Proc. Berkeley Conference in Honor of Jerzy Neyman and Jack Kiefer* (L. M. Le Cam and R. A. Olshen, eds.) 539–555. Wadswordth, Belmont, CA. MR0822052
[22] GROENEBOOM, P. (1989). Brownian motion with a parabolic drift and Airy functions. *Probab. Theory Related Fields* **81** 79–109. MR0981568
[23] GROENEBOOM, P., JONGBLOED, G. and WELLNER, J. A. (2001). A canonical process for estimation of convex functions: The "invelope" of integrated Brownian motion $+t^4$. *Ann. Statist.* **29** 1620–1652. MR1891741
[24] GROENEBOOM, P., JONGBLOED, G. and WELLNER, J. A. (2001). Estimation of a convex function: Characterizations and asymptotic theory. *Ann. Statist.* **29** 1653–1698. MR1891742
[25] HEILER, S. and WILLERS, R. (1988). Asymptotic normality of $R$-estimates in the linear model. *Statistics* **19** 173–184. MR0945375
[26] HERRNDORF, N. (1984). A functional central limit theorem for weakly dependent sequences of random variables. *Ann. Probab.* **12** 141–153. MR0723735
[27] HOLM, S. and FRISÉN, M. (1985). Nonparametric regression with simple curve characteristics. Research Report 4, Dept. Statistics, Univ. Göteborg.
[28] HUANG, Y. and ZHANG, C.-H. (1994). Estimating a monotone density from censored observations. *Ann. Statist.* **22** 1256–1274. MR1311975
[29] IBRAGIMOV, I. A. and LINNIK, Y. V. (1971). *Independent and Stationary Sequences of Random Variables*. Wolters-Noordhoff, Groningen. MR0322926
[30] JONGBLOED, G. (1995). Three statistical inverse problems. Ph.D. dissertation, Technical Univ. Delft.
[31] KIM, J. and POLLARD, D. (1990). Cube root asymptotics. *Ann. Statist.* **18** 191–219. MR1041391
[32] KOMLÓS, J., MAJOR, P. and TUSNÁDY, G. (1976). An approximation of partial sums of independent RV's, and the sample DF. II. *Z. Wahrsch. Verw. Gebiete* **34** 33–58. MR0402883
[33] LEURGANS, S. (1982). Asymptotic distributions of slope-of-greatest-convex-minorant estimators. *Ann. Statist.* **10** 287–296. MR0642740
[34] MAMMEN, E. (1991). Estimating a smooth monotone regression function. *Ann. Statist.* **19** 724–740. MR1105841
[35] MAMMEN, E. (1991). Regression under qualitative smoothness assumptions. *Ann. Statist.* **19** 741–759. MR1105842
[36] MCLEISH, D. L. (1975). A maximal inequality and dependent strong laws. *Ann. Probab.* **3** 829–839. MR0400382
[37] PELIGRAD, M. (1985). An invariance principle for $\phi$-mixing sequences. *Ann. Probab.* **13** 1304–1313. MR0806227
[38] POLLARD, D. (1984). *Convergence of Stochastic Processes*. Springer, New York. MR0762984
[39] PRAKASA RAO, B. L. S. (1969). Estimation of a unimodal density. *Sankhyā Ser. A* **31** 23–36. MR0267677
[40] ROBERTSON, T., WRIGHT, F. T. and DYKSTRA, R. L. (1988). *Order Restricted Statistical Inference*. Wiley, Chichester. MR0961262
[41] ROCKAFELLAR, R. T. (1997). *Convex Analysis*. Princeton Univ. Press. MR1451876
[42] SHORACK, G. R. and WELLNER, J. A. (1986). *Empirical Processes with Applications to Statistics*. Wiley, New York. MR0838963

ORDER RESTRICTED INFERENCE 59
[43] TAQQU, M. S. (1975). Weak convergence to fractional Brownian motion and to the Rosenblatt process. *Z. Wahrsch. Verw. Gebiete* **31** 287–302. MR0400329
[44] TAQQU, M. S. (1979). Convergence of integrated processes of arbitrary Hermite rank. *Z. Wahrsch. Verw. Gebiete* **50** 53–83. MR0550123
[45] VAN ES, B. JONGBLOED, G. and VAN ZUIJLEN, M. (1998). Isotonic inverse estimators for nonparametric deconvolution. *Ann. Statist.* **26** 2395–2406. MR1700237
[46] WANG, Y. (1994). The limit distribution of the concave majorant of an empirical distribution function. *Statist. Probab. Lett.* **20** 81–84. MR1294808
[47] WRIGHT, F. T. (1981). The asymptotic behaviour of monotone regression estimates. *Ann. Statist.* **9** 443–448. MR0606630

DEPARTMENT OF MATHEMATICAL SCIENCES
CHALMERS UNIVERSITY OF TECHNOLOGY
AND GÖTEBORG UNIVERSITY
SE-412 96 GÖTEBORG
SWEDEN
E-MAIL: dragi@math.chalmers.se

DEPARTMENT OF MATHEMATICAL STATISTICS
CENTRE FOR MATHEMATICAL SCIENCES
LUND UNIVERSITY
BOX 118
SE-221 00 LUND
SWEDEN
E-MAIL: Ola.Hossjer@mathstat.lu.se